\documentclass{amsart}
\usepackage{amsmath, amssymb, amsthm, amscd, graphicx}

\theoremstyle{plain}
\newtheorem{prop}{Proposition}[section]
\newtheorem{coro}[prop]{Corollary}

\newtheorem{lemm}[prop]{Lemma}

\theoremstyle{definition}
\newtheorem{defi}[prop]{Definition}
\newtheorem*{propA}{Property A}
\newtheorem*{propC}{Property C}
\newtheorem{nota}[prop]{Notation}
\newtheorem{exam}[prop]{Example}
\newtheorem{rema}[prop]{Remark}

\numberwithin{equation}{section}

\def\Reff#1; #2; #3; #4; #5; #6; #7\par{%
\bibitem{#1} #2, {\it #3}, #4 {\bf #5} (#6) #7}

\def\Ref#1; #2; #3; #4\par{%
\bibitem{#1} #2, {\it #3}, #4}

\catcode`\¨=\active\def¨{~} 

\renewcommand{\a}{\mathtt{a}}
\newcommand{\A}{\mathtt{A}}

\renewcommand{\b}{\mathtt{b}}

\newcommand{\bbb}[3]{b_{#1,#2}(#3)}

\newcommand{\bx}{x}
\newcommand{\by}{y}

\renewcommand{\c}{\mathtt{c}}
\newcommand{\card}{\mathtt{\#}}
\newcommand{\cc}{c}
\newcommand{\ccc}[2]{\cc(\DD{#1}{#2})}

\newcommand{\CG}{\Gamma}
\newcommand{\cl}[1]{[#1]}
\newcommand{\Class}{C}

\newcommand{\D}{\Delta}
\newcommand{\Dd}[2]{\delta_{#1,#2}}
\newcommand{\Ddt}[2]{\widetilde\delta_{#1,#2}}
\newcommand{\DD}[2]{\D_{#1}^{#2}}
\newcommand{\dd}{d}
\newcommand{\DDD}[2]{\Div(\D_{#1}^{#2})}
\newcommand{\DDDD}[2]{(\Div(\D_{#1}^{#2}),\nobreak\sm\nobreak)}
\renewcommand{\div}{\prec}
\newcommand{\Div}{\mathrm{Div}}
\newcommand{\DDiv}[1]{(\Div(#1),\nobreak\sm\nobreak)}
\newcommand{\dive}{\preccurlyeq}

\newcommand{\pz}{z}

\newcommand{\ea}{a}
\newcommand{\eb}{b}

\newcommand{\etc}{{\it etc.}}
\newcommand{\ev}{\rho}
\newcommand{\expo}{e}

\newcommand{\flip}[1]{\phi_{#1}}
\newcommand{\first}[1]{(#1)_{\!_1}}

\newcommand{\g}{\gamma}
\let\ge=\geqslant
\newcommand{\G}{\Gamma}
\renewcommand{\gcd}{\mathrm{gcd}}

\newcommand{\GG}[2]{\CG(\DD{#1}{#2})}

\newcommand{\hh}[1]{h_{#1}}
\newcommand{\hhh}[3]{\hh{#1}(\DD{#2}{#3})}

\newcommand{\ie}{{\it i.e.}}
\let\ince=\subseteq
\newcommand{\ind}{i}

\newcommand{\inv}{^{-1}}

\newcommand{\ir}{k}

\newcommand{\kk}{k}

\newcommand{\last}[1]{(#1)_{\!_\infty}}
\let\le=\leqslant
\newcommand{\lcm}{\mathrm{lcm}}

\newcommand{\mm}{m}

\newcommand{\nn}{n}
\newcommand{\NN}{N}
\newcommand{\NNN}{N}

\newcommand{\opp}{\mathord{\cdot}}

\newcommand{\pp}{p}
\newcommand{\ppp}{p}
\newcommand{\pppp}{p}

\newcommand{\qq}{q}

\newcommand{\resp}{{\it resp.{~}}}
\newcommand{\rr}{r}
\newcommand{\RR}{\mathbb{R}}

\let\s=\sigma

\newcommand{\sm}{<}
\newcommand{\sme}{\le}
\renewcommand{\sp}{s}
\newcommand{\spp}{t}
\renewcommand{\ss}[1]{\sigma_{#1}}
\newcommand{\sss}[2]{\sigma_#1^{#2}}
\renewcommand{\SS}[2]{S_{#1}^{#2}}
\newcommand{\SSS}[2]{\underline{S}_{#1}^{#2}}
\newcommand{\sx}{x}
\newcommand{\Sym}[1]{\mathfrak{S}_{#1}}

\newcommand{\uu}{u}

\newcommand{\vv}{v}

\newcommand{\ww}{w}
\newcommand{\WW}[2]{\vec w_{#1}^{#2}}

\newcommand{\xx}{x}
\newcommand{\XX}{X}

\newcommand{\yy}{y}

\begin{document}

\author{Patrick DEHORNOY}
\address{Laboratoire de Math\'ematiques Nicolas
Oresme UMR 6139\\ Universit\'e de Caen,
14032~Caen, France}
\email{dehornoy@math.unicaen.fr}
\urladdr{//www.math.unicaen.fr/\textasciitilde dehornoy}

\title{Still another approach to the braid ordering}

\keywords{braid group; orderable group;
well-ordering; normal form; fundamental braid}

\subjclass{20F36, 05A05, 20F60}

\begin{abstract}
We develop a new approach to the linear ordering of the
braid group~$B_\nn$, based on investigating its
restriction to the set¨$\DDD\nn\dd$ of all divisors
of~$\DD\nn\dd$ in the monoid¨$B_\infty^+$, \ie, to
positive $n$-braids whose normal form has length at
most¨$\dd$. In the general case, we compute several
numerical parameters attached with the finite orders
$\DDDD\nn\dd$. In the case of $3$~strands, we moreover
give a complete description of the increasing
enumeration of~$\DDDD3\dd$. We deduce a new and
specially direct construction of the ordering
on~$B_3$, and a new proof of the result that its
restriction to~$B_3^+$ is a well-ordering of ordinal
type~$\omega^\omega$.
\end{abstract}

\maketitle

The general aim of this paper is to investigate the
connection between the Garside structure of Artin's
braid groups and their distinguished linear ordering
(sometimes called the Dehornoy ordering). This leads
to a new, alternative construction of the ordering.

Artin's braid groups¨$B_\nn$ are endowed with several
interesting combinatorial structures. One of them
stems from Garside's analysis
\cite{Gar} and is nowadays known as a Garside
structure¨\cite{Dgk, McC}. It describes $B_\nn$ as the group
of fractions of a monoid¨$B_\nn^+$ with a rich divisibility
theory. One of the outcomes of this theory is a unique
normal decomposition for every braid in¨$B_\nn$ in
terms of simple braids, which are the divisors of
Garside's fundamental braid~$\D_\nn$, a finite
family of¨$B_\nn^+$ in one-to-one correspondence
with the permutations of $\nn$~objects. One obtains a
natural graduation of the monoid¨$B_\nn^+$ by
considering the family~$\DDD\nn\dd$ of all divisors
of~$\DD\nn\dd$, which also are the elements
of~$B_n^+$ whose normal form has length at
most~$\dd$.

On the other hand, the braid groups are equipped with a
distinguished linear ordering, which is compatible with
multiplication on the left, and admits a simple
combinatorial characterization¨\cite{Dfb}: a
braid¨$\bx$ is smaller than another braid¨$\by$ if,
among all expressions of the quotient¨$\bx\inv\by$ in
terms of the standard generators¨$\ss i$, there exists
at least one expression in which the
generator¨$\ss\mm$ with maximal (or minimal)
index¨$\mm$ appears only positively, \ie,
$\ss\mm$ occurs, but $\ss\mm\inv$ does not. Several
deep results about that ordering are known, in
particular the fact that its restriction
to¨$B_\infty^+$ is a well-ordering, and a number of
equivalent constructions are known¨\cite{Dgr}.

Although both combinatorial in nature, the previous
structures remain mostly uncon\-nected---and connecting
them may appear as one of the most natural questions of
braid combinatorics. For degree¨$1$, \ie, for simple
braids, the linear ordering corresponds to a
lexicographical ordering of the associated
permutations¨\cite{Dgb}. But this connection does not
extend to higher degrees, and almost nothing is known
about the restriction of the linear ordering to positive
braids of a given degree. In particular, no connection is
known between the above mentioned Garside normal form and
the alternative normal form constructed by S.\,Burckel
in¨\cite{Bus, But, Buu}, one that makes comparison with
respect to the linear ordering easy: to give an example,
the Garside normal form of¨$\DD3{2\dd}$
is¨$(\ss1\ss2\ss1)^{2\dd}$, while its Burckel normal form
is $(\ss2\ss1^2\ss2)^\dd\ss1^{2\dd}$.

Our aim in this paper is to investigate the finite
linearly ordered sets $\DDDD\nn\dd$. A
nice way of thinking of this structure is to consider the
increasing enumeration of¨$\DDD\nn\dd$, and to view it as
a distinguished path from¨$1$ to¨$\DD\nn\dd$ in the Cayley
graph of¨$B_\nn$. A complete description of this path
would arguably be an optimal solution to the rather vague
question of connecting the Garside and the ordered
structures of braid groups. Such a description seems to
be extremely intricate from a combinatorial point of
view, and it remains out of reach for the moment, but we
prove partial results in this direction, namely 

- $(i)$ in the general case, a determination of some
numerical parameters attached with $\DDDD\nn\dd$ that in
some sense measure its size, with explicit values for
small values of¨$\nn$ and¨$\dd$, and

- $(ii)$ in the special case $\nn = 3$, a complete
description of the increasing enumeration
of¨$\DDDD\nn\dd$.

More specifically, the parameters we investigate are
the complexity and the heights.
The complexity¨$\ccc\nn\dd$ is defined as the
maximal number of occurrences of¨$\ss{\nn-1}$ in an
expresion of¨$\DD\nn\dd$ containing no¨$\ss{\nn-1}$.
It is connected with the termination of the handle
reduction algorithm of¨\cite{Dfo}, and its
determination was left as an open question in the
latter paper. The
$\rr$-height¨$\hhh\rr\nn\dd$ is defined to
be the number of¨$\rr$-jumps in the increasing
enumeration of $\DDDD\nn\dd$ (augmented by¨$1$), where
the term $\rr$-jump refers to some natural filtration of
the linear ordering¨$\sm$ by a sequence of partial
orderings¨$\sm_\rr$. When $\rr$ increases, $\rr$-jumps
are higher and higher, so $\hhh\rr\nn\dd$ counts how
many big jumps exist in $\DDDD\nn\dd$.  We prove that the
complexity $\ccc\nn\dd$ equals the height
$\hhh{\nn-1}\nn\dd$ (Proposition¨\ref{P:MainHeight}), and
that, for each¨$\rr$, the $\rr$-height $\hhh\rr\nn\dd$ is
the number of divisors of¨$\DD\nn\dd$ whose $\dd$th
factor of the normal form is right divisible by¨$\D_\rr$
(Proposition¨\ref{P:Main}). Together with the
combinatorial results of¨\cite{Dhi}, this allows for
computing the explicit values listed in
Table¨\ref{T:Values}, and for establishing
various inductive formulas (Propositions¨\ref{P:Values}
and¨\ref{P:Values34}, among others).

Besides the enumerative results, we also prove a
general structural result that connects the ordered set
$\DDDD\nn\dd$ with (subsets of) $\DDDD{\nn-1}\dd$
(Corollary¨\ref{C:Structure}). This result
suggests an inductive method for directly constructing
the increasing enumeration of $\DDDD\nn\dd$ starting
from those of $\DDDD{\nn-1}\dd$ and $\DDDD\nn{\dd-1}$.
This approach is completed here for $\nn = 3$
(Proposition¨\ref{P:Enum3}). In some sense,
$3$~strand braids are simple objects, and the result
may appear as of modest interest; however, the order
on~$B_3^+$ is a well-ordering of ordinal
type~$\omega^\omega$, hence not a so simple object. The
interesting point is that this approach leads to a
new, alternative construction of the braid ordering,
with in particular a new and simple proof for the
so-called Comparison Property which is the hard core
in the construction, namely the part that guarantees
the linearity of the ordering. In this way, one
obtains not only one more construction of an ordering
that already has many constructions
\cite{Dgr}, but arguably the optimal one, as it makes
all proofs simple once the initial inductive
definition is correctly stated, and as the connection
with the Garside structure is then explicit.

The paper is organized as follows. After a first
introductory section recalling basic properties and setting
the notation, we introduce the parameters $\ccc\nn\dd$
and $\hhh\rr\nn\dd$ in Section¨\ref{S:Meas}, and we
establish their connection. In Section¨\ref{S:Deter}, we
connect in turn $\hhh\rr\nn\dd$ with the number of
$\nn$-braids whose $\dd$th factor in the normal form
satisfy certain constraints, and deduce explicit
values. Finally, in Section¨\ref{S:N3}, we study the
specific case of¨$\DDDD3\dd$ and describe its
increasing enumeration, resulting in the new
construction of the braid ordering in this case.

\section{Background and preliminary results}

Our notation is standard, and we refer to textbooks
like¨\cite{Bir} or¨\cite{Eps} for basic results about
braid groups. We recall that the $\nn$¨strand braid group¨$B_\nn$
is defined for $\nn \ge 1$ by the presentation
\begin{equation} \label{E:Present}
B_\nn = \left<\ss1, \dots, \ss{\nn-1} \,;\,
\begin{array}{cl}
\ss i \ss j = \ss j \ss i 
&\text{\quad for $\vert
i - j\vert \ge 2$}\\
\   \ss i \ss j \ss i  = \ss j \ss i \ss j
&\text{\quad for $\vert i - j\vert = 1$}
\end{array} \right>.
\end{equation}
So, $B_1$ is the trivial one-element group, while
$B_2$ is the free group generated by¨$\ss1$. The
elements of¨$B_\nn$ are called $\nn$¨strand braids,
or simply {\it $\nn$-braids}. We use $B_\infty$
for the group generated by an infinite sequence
of¨$\ss i$'s subject to the relations
of¨\eqref{E:Present}, \ie, the direct limit of
all¨$B_\nn$'s with respect to the inclusion of¨$B_\nn$
into¨$B_{n+1}$.

By definition, every $\nn$-braid¨$\bx$ admits (infinitely
many) expressions in terms of the generators¨$\ss i$, $1
\le i < n$. Such an expression is called an
$\nn$¨strand {\it braid word}. Two braid words¨$\ww,
\ww'$ representing the same braid are said to be {\it
equivalent}; the braid represented by a braid
word¨$\ww$ is denoted¨$\cl\ww$. 

\subsection{Positive braids and the element $\D_\nn$}

We denote by¨$B_\nn^+$ the monoid
admitting the presentation¨\eqref{E:Present},
and by¨$B_\infty^+$ the union (direct limit) of
all¨$B_\nn^+$'s. The elements of¨$B_\nn^+$ are called {\it
positive} $\nn$-braids. In¨$B_\infty^+$, no element
except¨$1$ is invertible, and we have a natural notion of
divisibility: 

\begin{defi}
For $\bx, \by$ in¨$B_\nn^+$, we say that $\bx$ is a {\it
left divisor} of¨$\by$, denoted $\bx \dive \by$, or,
equivalently, that $\by$ is a {\it right multiple}
of¨$\bx$, if $\by = \bx \pz$ holds for some¨$\pz$
in¨$B_\nn^+$. We denote by¨$\Div(\by)$ the (finite) set of
all left divisors of¨$\by$ in¨$B_\nn^+$.
\end{defi}

The monoid~$B_\nn^+$ is not commutative for $\nn \ge
3$, and therefore there are distinct symmetric
notions of a right divisor and a left multiple---but
we shall mostly use left divisors here. Note that
$\bx$ is a (left) divisor of¨$\by$ in the sense
of¨$B_\nn^+$ if and only if it is a (left) divisor in
the sense of¨$B_\infty^+$, so there is no need to
specify the index¨$\nn$.

According to Garside's theory \cite{Gar},
$B_\nn^+$ equipped with the left divisibility
relation is a lattice: any two
positive $\nn$-braids¨$\bx,
\by$ admit a greatest common left divisor, denoted
$\gcd(\bx, \by)$, and a least common right multiple,
denoted¨$\lcm(\bx, \by)$. A special role is played by the
lcm¨$\D_\nn$ of¨$\ss1$, \dots, $\ss{\nn-1}$, which can
be inductively defined by
\begin{equation} \label{E:Delta}
\D_1 = 1, \qquad
\D_\nn = \ss1 \ss2 \dots \ss{n-1} \, \D_{\nn-1}.
\end{equation}
It is well known that $\D_\nn^2$ belongs to the centre
of¨$B_\nn$ (and even generates it for $\nn \ge 3$),
and that the flip automorphism¨$\flip n$ of¨$B_\nn$
corresponding to conjugation by¨$\D_\nn$
exchanges¨$\ss i$ and¨$\ss{n-i}$ for
$1 \le i \le \nn-1$. It follows that $\flip n$
also exchanges¨$\ss{n, 1}$ and¨$\ss{1, n}$.

\begin{figure}[htb]
\begin{center}
\begin{picture}(83,21)(0,0)
\put(10,0){\includegraphics{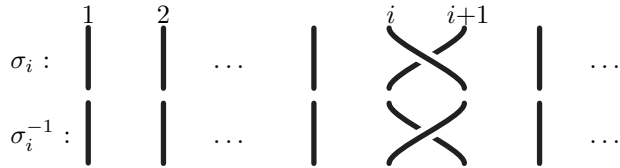}}
\put(0,13){$\ss i$~:}
\put(0,3){$\ss i\inv$~:}
\put(9.5,19){$1$}
\put(19.5,19){$2$}
\put(50,19){$i$}
\put(58,19){$i{+}1$}
\put(27,3){\dots}
\put(77,3){\dots}
\put(27,13){\dots}
\put(77,13){\dots}
\end{picture}
\end{center}
\caption{\smaller One associates with
every $\nn$¨strand braid word¨$\ww$ an $\nn$¨strand
braid diagram by stacking elementary diagrams as
above; then two braid words are equivalent if and only
if the associated diagrams are the projections of
ambient isotopic figures in¨$\RR^3$, \ie, one can
deform one diagram into the other without allowing the
strands to cross or moving the endpoints.}
\label{F:Diagram}
\end{figure}

In¨$B_\nn^+$, the left and the right divisors of¨$\D_\nn$
coincide, and they make a finite sublattice of¨$(B_\nn^+,
\dive)$ with $\nn!$¨elements. These braids will be called
{\it simple} in the sequel. When braid
words are represented by diagrams as
mentioned in Figure~\ref{F:Diagram}, simple braids
are those positive braids that can be represented by a
diagram in which any two strands cross at most
once.

By mapping¨$\ss i$ to the transposition¨$(i, i+1)$, one
defines a surjective homomorphism¨$\pi$ of¨$B_\nn$ onto the
symmetric group¨$\Sym\nn$. The restriction of¨$\pi$ to
simple braids is a bijection: for every permutation¨$f$
of¨$\{1, \dots, \nn\}$, there exists exactly one simple
braid¨$\bx$ satisfying $\pi(\bx) = f$. It follows that
the number of simple $\nn$-braids is¨$\nn!$. 

\begin{exam} \label{X:Degree}
The set¨$\DDD3{}$ consists of six elements, namely $1$,
$\ss1$, $\ss2$, $\ss2\ss1$, $\ss1\ss2$, and¨$\D_3$. In
examples, we shall often use the shorter notation¨$\tt a$
for¨$\ss1$, $\tt b$ for¨$\ss2$, \etc\ Thus, the six
simple $3$-braids are $1$, $\tt a$, $\tt b$, $\tt ba$,
$\tt ab$, $\tt aba$.  
\end{exam}

\subsection{The normal form}

For each positive $\nn$-braid¨$\bx$ distinct of¨$1$,
the simple braid $\gcd(\bx, \D_\nn)$ is the maximal
simple left divisor of¨$\bx$, and we obtain a
distinguished expression $\bx = \sx_1 \bx'$
with¨$\sx_1$ simple. By decomposing¨$\bx'$ in the
same way and iterating,  we obtains the so-called
normal expression¨\cite{ElM, Eps}. 

\begin{defi}
A sequence $(\sx_1, \dots, \sx_\dd)$ of simple
$\nn$-braids is said to be {\it normal} if, for
each¨$\kk$, one has $\sx_\kk =
\gcd(\D_\nn, \sx_\kk \dots \sx_\dd)$.
\end{defi}

Clearly, each positive braid admits a unique normal
expression. It will be convenient here to consider the
normal expression as unbounded on the right by
completing it with as many trivial factors¨$1$ as needed.
In this way, we can speak of the {\it $\dd$th factor} (in
the normal form) of¨$\bx$ for each positive braid¨$\bx$.
We say that a positive braid has {\it degree¨$\dd$} if
$\dd$ is the largest integer such that the $\dd$th factor
of¨$\bx$ is not¨$1$. We shall use the following two
properties of the normal form: 

\begin{lemm} \cite{ElM} \label{L:Normal}
A sequence of simple $\nn$-braids $(\sx_1, \dots,
\sx_\dd)$ is normal if and only if, for each¨$\kk <
\dd$, each¨$\ss i$ that divides¨$\sx_{\kk+1}$
on the left divides¨$\sx_\kk$ on the right.
\end{lemm}

\begin{lemm} \cite{ElM} \label{L:Degree}
For $\bx$ a positive braid in¨$B_\nn^+$, the following
are equivalent:

$(i)$ The braid¨$\bx$ belongs to¨$\DDD\nn\dd$, \ie, is a
(left or right) divisor of¨$\DD\nn\dd$; 

$(ii)$ The degree of¨$\bx$ is at most¨$\dd$.
\end{lemm}

\begin{exam} \label{X:Enum}
There are 19 divisors of~$\DD32$, which also are the
$3$-braids of degree at most~$2$. Their enumeration in
normal form---in an ordering that may seem strange now,
but should become familiar soon---is: $1$,
$\a$, $\a\opp\a$, $\b$, $\b\a$,
$\b\a\opp\a$, $\b\opp\b$, $\b\opp\b\a$, $\a\b$,
$\a\b\a$, $\a\b\a\opp\a$, $\a\b\opp\b$,
$\a\b\opp\b\a$, $\a\opp\a\b$, $\a\b\a\opp\b$,
$\a\b\a\opp\b\a$, $\b\a\opp\a\b$, $\a\b\a\opp\a\b$,
$\a\b\a\opp\a\b\a$. 
\end{exam}

By Lemma~\ref{L:Degree}, every divisor
of¨$\DD\nn\dd$ can be expressed as the product of at
most¨$\dd$ divisors of¨$\D_\nn$, so we certainly have
$\card\DDD\nn\dd \le (\nn!)^\dd$ for all¨$\nn,
\dd$.

\subsection{The braid ordering}

The basic notion is the following one:

\begin{defi}
Let $\ww$ be a nonempty braid word. We say that
$\ss\mm$ is the {\it main} generator in¨$\ww$
if $\ss\mm$ or $\ss\mm\inv$ occurs in~$\ww$,
but no $\ss i^{\pm1}$ with $i > \mm$ does. We
say that $\ww$ is {\it
$\s$-positive} (\resp {\it $\s$-negative}) if the main
generator occurs only positively (\resp negatively)
in¨$\ww$.
\end{defi}

A positive nonempty braid word, \ie, one that contains
no¨$\ss i\inv$ at all, is $\s$-positive, but the inclusion
is strict: for instance, $\ss1\inv \ss2$ is not positive,
but it is $\s$-positive, as its main generator,
namely¨$\ss2$, occurs positively (one¨$\ss2$) but not
negatively (no¨$\ss2\inv$).

Then we have the following two properties, which have
received a number of independent proofs¨\cite{Dgr}:

\begin{propA}
{\it A $\s$-positive braid word does not represent¨$1$.}
\end{propA}

\begin{propC}
{\it Every braid except¨$1$ can be represented by a
$\s$-positive word or by a $\s$-negative word.}
\end{propC}

It is then straightforward to order the braids:

\begin{defi}
If $\bx, \by$ are braids, we say that $\bx \sm \by$ holds
if the braid¨$\bx\inv\by$ admits at least one $\s$-positive
representative.
\end{defi} 

It is clear that the relation¨$\sm$ is transitive and
compatible with multiplication on the left;  Property¨A
implies that $\sm$ has no cycle, hence is a
strict partial order, and Property¨C then implies that it is
actually a linear order.

As every nonempty positive braid word is
$\s$-positive, $\bx \dive \by$ implies $\bx \sme \by$
for all positive braids¨$\bx, \by$, but the converse is
not true: $\ss1$ is not a left divisor of¨$\ss2$, but $\ss1
\sm \ss2$ holds, since $\ss1\inv\ss2$ is a $\s$-positive
word. 

\begin{exam} \label{X:List}
The increasing enumeration of the set¨$\DDD3{}$ is
$$\tt 1 \sm a \sm b \sm ba \sm ab \sm aba.$$
For
instance, we have $\tt ba \sm ab$, \ie, $\ss2\ss1
\sm \ss1\ss2$, as the quotient, namely $\ss1\inv
\ss2\inv \ss1 \ss2$ (or $\tt ABab$), also admits the
expression $\ss2 \ss1\inv$, a $\s$-positive word.
Similarly, the reader can check that the increasing
enumeration of¨$\DDD32$ is the one given in
Example~\ref{X:Enum}.
\end{exam}

\begin{lemm} \label{L:Extend}
The linear ordering¨$\sm$ extends the left divisibility
ordering¨$\div$.
\end{lemm}

\begin{proof}
By definition, $1 \sm \ss i$ holds for every¨$i$. As the
ordering¨$\sm$ is compatible with multiplication on the
left, it follows that $\xx \sm \xx \ss i$ holds for
all¨$i, \xx$, and, therefore, $\xx \sm \xx \yy$ holds
whenever $\yy$ is a non-trivial positive braid.
\end{proof}

Lemma¨\ref{L:Extend} implies that $1$ is always the
first element of¨$\DDDD\nn\dd$, and
$\DD\nn\dd$ is always its last element.
It may be noted that a deep result by
Laver¨\cite{Lve} shows that, although $\sm$ is not
compatible with right multiplication in general,
nevertheless $\xx \sm \ss i
\xx$ always holds,
\ie, $\sm$ also extends the right divisibility ordering. 

By Property~C, every nontrivial braid admits
at least one $\s$-positive or $\s$-negative
expression. In general, such a $\s$-positive or
$\s$-negative expression is not unique, but the main
generator in such expressions is uniquely defined:

\begin{lemm} \label{L:Index}
If a braid¨$\bx$ admits a $\s$-positive 
expression, then the main generators in any
two $\s$-positive expressions of¨$\bx$
coincide.
\end{lemm}

\begin{proof}
Assume that $\ww$, $\ww'$ are
$\s$-positive  expressions of¨$\bx$, and let
$\ss\mm$, $\ss{\mm'}$ be their main
generators. Assume for instance $\mm <
\mm'$. Then $\ww\inv \ww'$ is a
$\s$-positive word, and it represents the
trivial braid¨$1$: this contradicts Property¨A.
\end{proof}

So, there will be no ambiguity in referring to
{\it the} main generator of some non-trivial braid¨$\bx$: 
this means the main generator in any $\s$-positive (or
$\s$-negative) expression of¨$\bx$.

\begin{rema}
Our current definition corresponds to the order¨$<^\phi$
of¨\cite{Dgr}. It differs from the one most usually
considered in literature in that we refer to the maximal
index rather than to the minimal one in the
definition of a $\s$-positive word. Switching from one
definition to the other amounts to conjugating
by¨$\D_\nn$, \ie, to applying the flip automorphism. Results
are entirely similar for both versions. However, it is much
more convenient to consider the ``max'' choice here, because
it guarantees that $B_\nn^+$ is an initial segment
of¨$B_{n+1}^+$. Using the ``min'' convention would make the
statements of the forthcoming sections less natural.
\end{rema}

\section{Measuring the ordered sets¨$\DDDD\nn\dd$}
\label{S:Meas}

As was said above, our aim is
to investigate the finite ordered sets¨$\DDDD\nn\dd$,
and, more generally, $(\Div(\pz), \sm)$ for¨$\pz$ a
positive braid. We shall do it by defining numerical
parameters that somehow measure their size. The first
parameter involves the length of the
$\s$-positive words that are, in some natural sense
defined below, drawn in the Cayley graph
of¨$\DD\nn\dd$. It will be called the {\it complexity}
of¨$\DD\nn\dd$, because it is directly connected with
the complexity analysis of the handle reduction
algorithm of¨\cite{Dfo}. The other parameters involve
a filtration of the linear ordering by
the¨$\ss i$'s, and they will be called the {\it
heights} of¨$\DD\nn\dd$ because they count the jumps
of a given height in¨$\DDDD\nn\dd$. 

\subsection{Sigma-positive paths in the Cayley graph}

The first parameter we attach to $(\Div(\pz),
\sm)$ involves the $\s$-positive paths in the Cayley
graph of¨$\pz$.

We recall that the Cayley graph of
the group¨$B_\nn$ with respect to the standard
generators $\ss i$ is the labeled graph  with vertex
set¨$B_\nn$ and such that there exists a $\ss
i$-labeled edge from¨$\bx$ to¨$\by$ if and only if
$\by = \bx \ss i$ holds. The Cayley graph of the
monoid¨$B_\nn^+$ is obtained by restricting the vertices
to¨$B_\nn^+$. Note that the Cayley graph of¨$B_\nn$
(and a fortiori of¨$B_\nn^+$) can be seen as a subgraph
of the Cayley graph of¨$B_\infty$.

\begin{figure} [htb]
\begin{picture}(110,51)(0, -3)
\put(5,8){\includegraphics{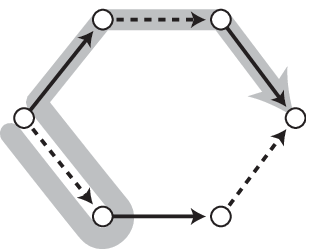}}
\put(50,0){\includegraphics{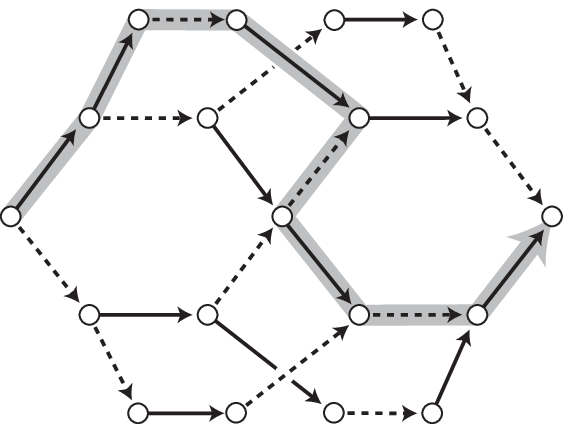}}
\put(2,20){$1$}
\put(10,34){$\ss2$}
\put(27,34){$\ss2\ss1$}
\put(10,7){$\ss1$}
\put(27,7){$\ss1\ss2$}
\put(36,20){$\D_3$}
\put(47,20){$1$}
\put(55,33){$\ss2$}
\put(60,44){$\ss2^2$}
\put(70,44){$\ss2^2\ss1$}
\put(80,44){$\ss2\ss1^2$}
\put(90,44){$\ss2\ss1^2\ss2$}
\put(65,33){$\ss2\ss1$}
\put(65,33){$\ss2\ss1$}
\put(86,33){$\D_3\ss1$}
\put(99,33){$\D_3\ss1\ss2$}
\put(55,8){$\ss1$}
\put(65,8){$\ss1\ss2$}
\put(86,8){$\D_3\ss2$}
\put(99,8){$\D_3\ss2\ss1$}
\put(82,20){$\D_3$}
\put(108,20){$\D_3^2$}
\put(60,-3){$\ss1^2$}
\put(70,-3){$\ss1^2\ss2$}
\put(80,-3){$\ss1\ss2^2$}
\put(90,-3){$\ss1\ss2^2\ss1$}
\end{picture}
\caption{\smaller The graphs
of¨$\GG3{}$ and $\GG32$;
the dotted edges represent¨$\ss1$, the plain
ones¨$\ss2$; observe that the graph
of¨$\DD32$ is not planar; in grey: two
$\s$-positive words traced in the graphs, namely
$\tt aAbab$ and $\tt bbabAbab$ ({\it cf.}
Lemma¨\ref{L:Path})}
\label{F:Cayley}
\end{figure}

\begin{defi}
(Figure¨\ref{F:Cayley})
For $\pz$ a positive braid, we denoted by¨$\CG(\pz)$
the subgraph of the Cayley graph of¨$B_\infty$
obtained by restricting the vertices to¨$\Div(\pz)$,
and only keeping those edges that connect two vertices
in¨$\Div(\pz)$.
\end{defi}

As every element of¨$B_\nn^+$ is a left divisor of¨$\DD\nn\dd$
for¨$\dd$ large enough, the Cayley graph of¨$B_\nn^+$ is 
the union of all graphs¨$\GG n k$ when $\dd$¨varies.

A path in the Cayley graph can be specified by its
initial vertex and the list of the labels of its
successive edges, \ie, by a braid word. For each¨$i <
n$ and each¨$\bx$ in¨$B_\nn$, there is exactly one $\ss
i$-labeled edge with target¨$\bx$, and one
$\ss i$-labeled edge with source¨$\bx$ in the Cayley
graph of¨$B_\nn$. Hence, in the complete
Cayley graph of¨$B_\nn$, for each initial vertex¨$\bx$ and
each $\nn$-braid word¨$\ww$, there is always one path
labeled¨$\ww$ starting from¨$\bx$. When we restrict to some
fragment¨$\G$, this need not be the case, but we have an
unambiguous notion of¨$\ww$ being drawn in¨$\G$
from¨$\bx$. Formally, this corresponds to

\begin{defi}
For $\G$ a subgraph of the Cayley graph
of¨$B_\infty$, and $\bx$ a vertex in¨$\G$, we say
that a braid word¨$\ww$ is {\it drawn} from¨$\bx$
in¨$\G$ if, for every prefix¨$\vv\ss i$ (\resp $\vv\ss
i\inv$) of¨$\ww$, there exists a $\ss i$-labeled edge
starting from (\resp finishing at)¨$\bx
\, \cl{\vv}$ in¨$\G$.
\end{defi}

For instance, we can check on Figure¨\ref{F:Cayley} that
the word¨$\ss1^2$ is drawn from¨$\ss2$
in¨$\GG32$, but not in¨$\GG3{}$. In algebraic terms,
we have the following characterization:

\begin{lemm} \label{L:Path}  
Assume that $\pz$ is a positive braid, and $\ww$ is a
braid word. Then $\ww$ is drawn from¨$\bx$ in¨$\CG(\pz)$
if and only if $\bx \cl{\vv} \dive \pz$ holds for each
prefix¨$\vv$ of¨$\ww$.
\end{lemm}               

\begin{proof}
The condition is sufficient. Indeed, assume it is
satisfied by¨$\ww$, and $\vv \ss i$ is a prefix of¨$\ww$.
Then, by hypothesis, $\bx \cl{\vv}$ and
$\bx \cl{\vv}\ss i$ are left divisors of¨$\pz$, hence
are vertices in¨$\CG(\pz)$, and, therefore,
there is a $\ss i$-labeled edge between¨$\bx \cl{\vv}$ and
$\bx \cl{\vv} \ss i$ in¨$\CG(\pz)$. The argument is
similar for a prefix of the form $\vv \ss i\inv$. Using
induction on the length of¨$\ww$, we deduce that $\ww$ is
drawn from¨$\bx$ in¨$\CG(\pz)$. 

Conversely, if there exists a $\ww$-labeled path
from¨$\bx$ in¨$\CG(\pz)$, then, for each prefix¨$\vv$
of¨$\ww$, the braid¨$\bx \cl{\vv}$ has to represent some
vertex in¨$\CG(\pz)$, hence it is a left divisor
of¨$\pz$. 
\end{proof}

For¨$\pz$ a positive braid, we shall investigate the
$\s$-positive words drawn in the graph¨$\CG(\pz)$. It is
clear that, even if $\Div(\pz)$ is a finite set,
arbitrary long words are drawn in¨$\CG(\pz)$ whenever the
latter contains at least 2¨vertices, \ie, $\pz$ is
not¨$1$. The example of Figure¨\ref{F:Cayley} shows that
restricting to $\s$-positive words does not change the
result: for instance, for each¨$\kk$, the
word $(\ss1\ss1\inv)^\kk\ss2\ss1\ss2$ is a
$\s$-positive expression of¨$\D_3$, and it is drawn
in¨$\GG3{}$. So we cannot hope for any finite upper
bound on the length of the $\s$-positive words drawn
in¨$\CG(\pz)$ in general. Now, the situation changes
if we concentrate on the main generators, \ie, we
forget about the generators with non-maximal index.

\begin{lemm} \label{L:Finite}
Assume that $\G$ is subgraph of the Cayley graph
of¨$B_\infty$, and $\ww$ is a $\s$-positive word drawn
in¨$\CG(\pz)$. Then the number of occurrences of the main
generator in¨$\ww$ is at most the number of
non-terminal vertices in¨$\G$.
\end{lemm}

\begin{proof}
Assume that $\ww$ is drawn from¨$\bx$ in¨$\G$. Let
$\ss\mm$ be the main generator in¨$\ww$. As there is at
most one $\ss\mm$-labeled edge starting from each vertex
of¨$\G$, it suffices to show that the number of¨$\ss\mm$'s
in¨$\ww$ is bounded above by the number of
$\ss\mm$-edges in¨$\G$. Hence, it suffices to show that
the path¨$\gamma$ associated with¨$\ww$ cannot cross the
same $\ss\mm$-edge twice. Now assume that some
$\ss\mm$-edge starts from the vertex¨$\by$, and that $\g$
crosses this edge twice. This means that $\g$ contains a
loop from¨$\by$ to¨$\by$. Let
$\vv$ be the subword of¨$\ww$ labeling that loop. By
construction, $\vv$ begins with¨$\ss\mm$, it contains
no¨$\ss\mm\inv$ and no $\ss i^{\pm1}$ with $i > m$ as
it is a subword of¨$\ww$, and it represents the
braid¨$1$ as it labels a loop in the Cayley graph
of¨$B_\infty$: this means that $\vv$ is a
$\s$-positive word representing¨$1$, which contradicts
Property¨A.
\end{proof}

Lemma¨\ref{L:Finite} applies in particular to every
graph¨$\CG(\pz)$ with¨$\pz$ a positive braid. So we can
introduce our first parameter measuring
the size of the ordered set¨$\DDiv\pz$:

\begin{defi}
(Figure~\ref{F:Cayley})
For¨$\pz$ a positive braid with main generator¨$\ss\mm$,
the {\it complexity}¨$\cc(\pz)$ of¨$\pz$ is defined to
be the maximal number of¨$\ss\mm$'s in a $\s$-positive
word drawn in¨$\CG(\pz)$. 
\end{defi}

\begin{exam}
The word¨$\ss2 \ss1 \ss2$ is a
$\s$-positive word drawn from¨$1$ in¨$\GG3{}$, and
it contains two¨$\ss2$'s, hence we have $\ccc3{}
\ge 2$. Actually, it is not hard to obtain the
exact value $\ccc3{} = 2$. Indeed, if a
$\s$-positive path¨$\g$ contains the two
$\ss2$-edges starting from¨$1$ and¨$\ss1\ss2$, it
cannot come back to¨$\ss2$ for possibly crossing the
third $\ss2$-edge; and if $\g$ contains the
$\ss2$-edge that starts from¨$\ss1$, it can never
come back to¨$1$ or to¨$\ss2\ss1$ and therefore
contains at most one $\ss2$-edge. 
As we have $\DD3\dd = (\ss2 \ss1 \ss2)^\dd$,
we deduce $\ccc3\dd \ge 2\dd$ for every¨$\dd$;
this value is certainly not optimal, since the
example displayed in Figure¨\ref{F:Cayley} contains
five¨$\ss2$'s, proving $\ccc32 \ge 5$---the exact
value is¨$6$, and, more generally, we have $\ccc\nn\dd
= 2^{\dd+1} - 2$, as will be seen in
Section¨\ref{S:Deter}.
\end{exam}

\begin{rema}
Restricting to $\s$-positive words drawn in¨$\CG(\pz)$
is essential: for instance, for each¨$\kk$, we have
\begin{equation} \label{E:Unbounded}
\D_3 = \ss2^{\kk+1} \ss1\ss2\ss1^{-\kk},
\end{equation}
a $\s$-positive word containing $\kk + 2$ letters¨$\ss2$.
Now, for¨$\kk \ge 1$, the word involved
in¨\eqref{E:Unbounded} is not drawn in¨$\GG3{1}$, as
its prefix¨$\ss2^2$ is not. Thus the
parameter¨$\cc(\pz)$ really involves the left
divisors of¨$\pz$.
\end{rema}

A direct application of Lemma¨\ref{L:Finite} gives:

\begin{prop}
Every positive braid has a finite complexity; more
precisely, for¨$\pz$ of length¨$\ell$ in¨$B_\nn^+$  with
$\nn \ge 3$, we have $\cc(\pz) \le (\nn-1)^\ell$.
\end{prop}

\begin{proof}
The number of non-terminal vertices in¨$\CG(\pz)$, \ie,
the number of proper left divisors of¨$\pz$, is at most
$1 + (\nn-1) + (\nn-1)^2 + \cdots + (\nn-1)^{\ell-1}$.
\end{proof}

As the length of any positive expression of¨$\D_\nn$ is
$\nn(\nn-1)/2$, we obtain in particular for all¨$\nn, \dd$
\begin{equation} \label{E:Coarse}
\ccc\nn\dd \le
(\nn-1)^{\dd\nn(\nn-1)/2}.
\end{equation}

Before going further, let us observe that, in the
definition of the complexity of¨$\pz$, we can restrict
to decompositions of¨$\pz$, \ie, instead of considering
paths starting from and finishing at arbitrary vertices,
we can restrict to paths starting from¨$1$ and finish
at¨$\pz$:

\begin{lemm} \label{L:Decomp}
Assume that $\pz$ is a positive braid with main
generator¨$\ss\mm$. Then $\cc(\pz)$ is the
maximal number of¨$\ss\mm$'s in any $\s$-positive
decomposition of¨$\pz$ drawn in¨$\CG(\pz)$.
\end{lemm}

\begin{proof}
Let $\cc'(\pz)$ be the number involved in the above
statement. Clearly we have $\cc'(\pz) \le \cc(\pz)$.
Conversely, assume that $\ww$ is drawn in¨$\CG(\pz)$
from¨$\bx$, and that the $\ww$-labeled path starting
from¨$\bx$ finishes at¨$\by$. Let $\uu$ be a positive
expression of¨$\bx$, and $\vv$ be a positive expression
of¨$\by\inv \pz$. The latter exists as, by hypothesis,
$\by$ is a left divisor of¨$\pz$. Then $\uu\ww\vv$ is a
$\s$-positive decomposition of¨$\pz$ drawn in¨$\CG(\pz)$.
Hence we have $\cc'(\pz) \ge \cc(\pz)$.
\end{proof}

\begin{rema}
Let us call Property~A$^*$ the fact that all
numbers¨$\ccc\nn\dd$ are finite. Above we derived 
Property~A$^*$ from Property¨A. Actually, the
implication is an equivalence, \ie, we can also
deduce Property¨A from Property~A$^*$. Indeed,
assume that some $\s$-positive braid word¨$\ww$
represents¨$1$. The word¨$\ww$ may involve negative
letters and the problem is to find a vertex¨$\bx$
such that there exists a path labeled¨$\ww$
from¨$\bx$ in some¨$\GG\nn\dd$. Let $\ss m$ be the
main generator in¨$\ww$. The
word¨$\ww$ has finitely many prefixes, say
$\ww_0$, \dots, $\ww_\ell$. By Garside's theory, each
word¨$\ww_i$ is equivalent to a word of the form¨$u_i\inv
\vv_i$ with $\uu_i, \vv_i$ positive. Let $\bx$
be the least common left multiple of the
positive braids¨$\cl{\uu_0}$,
\dots, $\cl{\uu_\ell}$. Then, for each¨$i$, the
braid¨$\bx\cl{\ww_i}$ is positive. Moreover, there
exist¨$\nn$ and¨$\dd$ such that $\bx\cl{\ww_0}$, \dots,
$\bx\cl{\ww_\ell}$ all are divisors of¨$\DD\nn\dd$. This
means that the word¨$\ww$ is drawn from¨$\bx$ 
in¨$\GG\nn\dd$, and the associated path is a loop
around¨$\bx$. It follows that $\ww^\kk$ is drawn
in¨$\GG\nn\dd$ from¨$\bx$ for each¨$\kk$. By
construction, $\ww^\kk$ contains at least
$\kk$¨generators¨$\ss m$, hence $\ccc\nn\dd$ cannot
be finite.
\end{rema}

\subsection{Connection with handle reduction}

Handle reduction¨\cite{Dfo} is an algorithmic solution to
the word problem of braids that relies on the braid
ordering---actually the most efficient method
available to-date in practice. It is proved to be
convergent, but the complexity upper bound resulting
from the argument of¨\cite{Dfo} is exponential with
respect to the length of the input word, seemingly
very far from sharp. 

Each step of handle reduction involves a specific
generator¨$\ss i$, and, for an induction, the point is to
obtain an upper bound on the number of reduction steps
involving the main generator. The latter will naturally
be called the {\it main} reduction steps. The
connection between handle reduction and the complexity
as defined above relies on the following
technical result: 

\begin{lemm} \cite{Dfo}
Assume that $\pz$ is a positive braid with main
generator¨$\ss\mm$, and $\ww$ is drawn in¨$\CG(\pz)$.
Then, for each sequence of handle reductions from¨$\ww$,
\ie, each sequence¨$\vec\ww$ with $\ww_0= \ww$ such that
$\ww_\kk$ is obtained by reducing one handle
from¨$\ww_{\kk-1}$ for each¨$\kk$, there exists a
witness-word¨$\uu$ that is  $\s$-positive, drawn
in¨$\Div(\pz)$, and such that the number of¨$\ss\mm$'s
in¨$\uu$ is the number of main reductions in¨$\vec\ww$.
\end{lemm}

It follows that the number of main reduction steps in any
sequence of handle reductions starting with a word drawn
in¨$\CG(\pz)$ is bounded above by¨$\cc(\pz)$. In
particular, if we start with an $n$¨strand braid
word¨$\ww$ of length¨$\ell$, then it is easy to show
that $\ww$ is drawn in¨$\GG\nn\ell$, and, applying the
upper bound of¨\eqref{E:Coarse}, we deduce an
upper bound for the number of possible main reductions
from¨$\ww$, one that is exponential with respect
to¨$\ell$.

A natural way of improving this coarse upper bound would
be to determine the value of¨$\ccc\nn\dd$ more precisely.
This will be done in Section¨\ref{S:Deter} below.
However, the explicit formulas show that, for $\nn \ge
3$, the growth rate with respect to¨$\dd$ is exponential,
thus discarding any hope of proving the expected
polynomial upper bound for the number of reduction steps
by this approach.

\subsection{A filtration of the braid ordering}

We now introduce new numerical parameters for the ordered
sets¨$\DDiv\pz$. These numbers appear in connection with a
natural filtration of the ordering¨$\sm$, using an
increasing sequence of partial orderings that we
introduce now.

By Lemma¨\ref{L:Index}, the index of the main generator
of a non-trivial braid is well defined. We can use this
index to  measure the height of the jump between two
braids¨$\bx, \by$ satisfying $\bx \sm \by$:

\begin{defi}
For $\bx, \by$ in~$B_\infty$ and $\rr
\ge 1$, we say that $\bx \sm_\rr \by$ holds
or, equivalently, that $(\bx, \by)$ is a
{\it $\rr$-jump}, if $\bx\inv\by$ admits a
$\s$-positive expression in which the main
generator is~$\ss\mm$ with $\mm \ge \rr$.
\end{defi}

\begin{lemm}
For each $\rr \ge 1$, the
relation~$\sm_\rr$ is a strict partial
order that refines~$\sm$; the
relation~$\sm_1$ coincides with~$\sm$,
and $\rr \le \qq$ implies that $\sm_\qq$
refines~$\sm_\rr$. 
\end{lemm}

\begin{proof}
That $\sm_\rr$ is transitive follows from
the fact that the concatenation of a
$\s$-positive word with main
generator~$\ss{\mm}$ and a $\s$-positive word 
with main generator~$\ss{\mm'}$ is a
$\s$-positive word with main
generator~$\ss{\max(\mm, \mm')}$. 
\end{proof}

In the sequel, we consider the $\sm_\rr$-chains
included in~$\Div(\pz)$, and their length:

\begin{defi}
For $\pz$ a positive braid and $\rr \ge 1$,
we define the {\it $\rr$-height}¨$\hh\rr(\pz)$ of¨$\pz$
to be the maximal length of a $\sm_\rr$-chain included
in~$\Div(\pz)$. 
\end{defi}

Before giving examples, we observe the
connection between¨$\hh\rr(\pz)$ and
the increasing enumeration of the
set~$\Div(\pz)$:

\begin{lemm} \label{L:Enum}
Let $\pz$ be a positive braid and $\rr \ge
1$. Then $\hh\rr(\pz)-1$ is the
number of $\rr$-jumps in the increasing
enumeration of~$\DDiv\pz$.
\end{lemm}

\begin{proof}
If the number of $\rr$-jumps in the
increasing enumeration of~$\Div(\pz)$
is~$\NN_\rr-1$, we can extract from~$\Div(\pz)$ a
$\sm_\rr$-chain of length~$\NN_\rr$. Conversely,
assume that $(\by_0, \dots \by_{\NN_\rr})$ is a
$\sm_\rr$-chain in~$\Div(\pz)$. 
Let $\pz_0 \sm \ldots \sm \pz_\NNN$ be the
increasing enumeration of~$\Div(\pz)$. As
$\sm_\rr$ refines~$\sm$, there exists an
increasing function~$f$ of $\{0, \ldots,
\NN_\rr\}$ into $\{0, \ldots, \NNN\}$ such that
$\by_i = \pz_{f(i)}$ holds for every~$i$. Now
the hypothesis $\pz_{f(i)} \sm_\rr
\pz_{f(i+1)}$ implies that there exists at
least one $\rr$-jump between $\pz_{f(i)}$
and~$\pz_{f(i+1)}$. Indeed, by
Lemma~\ref{L:Index}, it is impossible that a
concatenation of $\mm$-jumps with $\mm < \rr$
results in a $\rr$-jump. So the number
of $\rr$-jumps in $(\pz_0, \dots, \pz_\NNN)$ is
at least~$\NN_\rr$.
\end{proof}

In other words, in order to determine¨$\hh\rr(\pz)$,
there is no need to consider arbitrary chains: it
is enough to consider the maximal chain
obtained by enumerating $\Div(\pz)$ exhaustively. 

\begin{exam} \label{X:Height}
Refining the increasing enumeration of¨$\DDD3{}$
given in Example¨\ref{X:List} by indicating for
each step the height of the corresponding jump, we
obtain:
\begin{equation} \label{E:List31}
\tt 1 \sm_1 a \sm_2 b
\sm_1 ba
\sm_2 ab \sm_1 \D_3,
\end{equation}
where we recall $\tt a, b, \dots$ stand for
$\ss1, \ss2, \dots$. For instance,
$\tt (ba, ab)$ is a $2$-jump, as we have $\tt
(ba)\inv(ab) = ABab = AabA = bA$, a
$\s$-positive decomposition with main
generator~$\ss2$. The number of $1$-jumps
in~\eqref{E:List31},
\ie, the number of symbols~$\sm_\rr$ with
$\rr \ge 1$, is~$5$, while
the number of~$2$-jumps is~$2$, so, by
Lemma¨\ref{L:Enum}, we deduce
$\hhh13{} = 6$ and $\hhh23{} = 3$.
Similarly, we obtain for¨$\DD32$ 
\begin{multline*}
\tt 1 \sm_1 a \sm_1 aa \sm_2 b \sm_1 ba 
\sm_1 baa \sm_2 bb \sm_1 bba 
\sm_2 ab \sm_1
aba  \sm_1 abaa \sm_2 abb 
\\ \tt 
\sm_1 abba \sm_2
aab
\sm_1  aaba \sm_1 aabaa \sm_2 baab \sm_1 baaba
\sm_1 baabaa,
\end{multline*}
leading to $\hhh132 = 19$ and $\hhh232 = 7$.
\end{exam}

\begin{prop} \label{P:Decreasing}
$(i)$ For every braid~$\pz$ in~$B_\nn^+$, we
have
\begin{equation}
\hh1(\pz) = \card\Div(\pz)  \ge
\hh2(\pz)
\ge \dots \ge \hh\nn(\pz) = 1.
\end{equation}

$(ii)$ For all positive braids~$\pz, \pz'$
and $\rr \ge 1$, we have
\begin{equation} \label{E:Addition}
\hh\rr(\pz\pz') \ge \hh\rr(\pz) + \hh\rr(\pz').
\end{equation}
\end{prop}

\begin{proof}
$(i)$ A $\sm_1$-chain is simply a $\sm$-chain,
hence every subset of~$\Div(\pz)$ gives such
a chain. So the maximal $\sm_1$-chain in
$\Div(\pz)$ is $\Div(\pz)$ itself, and
$\hh1(\pz)$ is the cardinality of~$\Div(\pz)$.

On the other hand, no $\sm_\nn$-chain
in~$B_\nn^+$ has length more than~$1$, as the
main generator of a $\s$-positive $\nn$~strand
braid word cannot be $\ss\nn$ or above. So
$\hh\nn(\pz)$ is~$1$.

Then, for $\qq \le \rr$, every
$\sm_\rr$-chain is a $\sm_\qq$-chain, which
implies $\hh\rr(\pz) \ge \hh\qq(\pz)$.

Point~$(ii)$ is obvious, as the
concatenation of two $\sm_\rr$-chains is a
$\sm_\rr$-chain.
\end{proof}

From~\eqref{E:Addition} we deduce
$\hh\rr(\pz^\dd) \ge \dd \cdot \hh\rr(\pz)$
for all~$\rr, \pz$. By Lemma¨\ref{L:Degree}, every
divisor of~$\DD\nn\dd$ can be decomposed as the
product of at most~$\dd$ divisors
of~$\D_\nn$, and the latter are
$\nn!$~in number, so we obtain the (coarse) bounds
\begin{equation}
\dd \cdot \hhh\rr\nn{} \le \hhh\rr\nn\dd \le
(\nn!)^\dd
\end{equation}
for all~$\rr, \nn, \dd$. Better estimates
will be given below.

\begin{rema}
Instead of restricting to subsets of~$B_\infty$ of the
form~$\Div(\pz)$, we can define the complexity and the
$\rr$-height for every (finite) set of braids~$\XX$. Most
general results extend, but, when $\XX$ is not closed
under left division, nothing can be said about the
number of~$\ss\rr$'s involved in a
$\rr$-jump. Considering such an extension is not useful
here. 
\end{rema}

\subsection{Connection with the complexity}

We shall now connect the complexity¨$\cc(\pz)$ with
the numbers¨$\hh\rr(\pz)$ just defined. The result
is simple:

\begin{prop} \label{P:MainHeight}
For $\pz$ a positive braid with main
generator¨$\ss\mm$, we have
\begin{equation}
\cc(\pz) = \hh\mm(\pz) - 1.
\end{equation}
In particular, for $\nn \ge 2$ and $\dd
\ge 0$, we have 
\begin{equation} \label{E:Complexity}
\ccc\nn\dd = \hhh{\nn-1}\nn\dd - 1.
\end{equation}
\end{prop}

One inequality is easy:

\begin{lemm} \label{L:MainHeight}
For $\pz$ a positive braid with main
generator¨$\ss\mm$, we have $\cc(\pz) \le
\hh\mm(\pz) - 1$.
\end{lemm}

\begin{figure} [htb]
\begin{picture}(120,27)(0, -0)
\put(0,2){\includegraphics{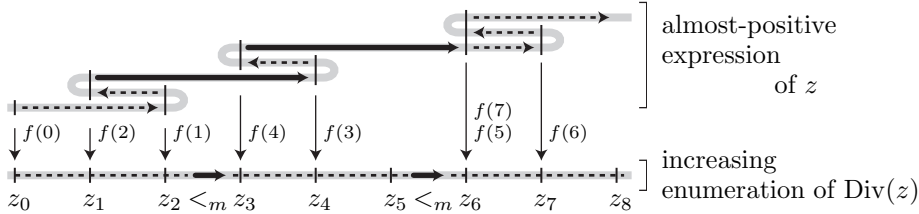}}
\put(87,23){almost-positive}
\put(87,19){expression}
\put(102,15){of¨$\pz$}
\put(87,5){increasing }
\put(87,1){enumeration of $\Div(\pz)$}
\put(0,0){$\pz_0$}
\put(10,0){$\pz_1$}
\put(20,0){$\pz_2$}
\put(24,0){$\sm_\mm$}
\put(30,0){$\pz_3$}
\put(40,0){$\pz_4$}
\put(50,0){$\pz_5$}
\put(54,0){$\sm_\mm$}
\put(60,0){$\pz_6$}
\put(70,0){$\pz_7$}
\put(80,0){$\pz_8$}
\put(2,9){$\scriptstyle f(0)$}
\put(12,9){$\scriptstyle f(2)$}
\put(22,9){$\scriptstyle f(1)$}
\put(32,9){$\scriptstyle f(4)$}
\put(42,9){$\scriptstyle f(3)$}
\put(62,9){$\scriptstyle f(5)$}
\put(62,12){$\scriptstyle f(7)$}
\put(72,9){$\scriptstyle f(6)$}
\end{picture}
\caption{\smaller Proof of
Lemma~\ref{L:MainHeight}: the main
generator~$\ss\mm$ corresponds to the bold
arrow: the function¨$f$ need not be
increasing, but the projection of a bold
arrow upstairs must include at least one bold
arrow downstairs, \ie, at least one
$\mm$-jump.}
\label{F:Jump}
\end{figure}

\begin{proof} 
The argument is reminiscent of that used for
Lemma~\ref{L:Enum}, but requires a little more care.
Assume that $\ww$ is a $\s$-positive word drawn
in¨$\CG(\pz)$ from¨$\bx$ containing
$\NN_\mm$¨occurrences of¨$\ss\mm$. By
Lemma¨\ref{L:Decomp}, we can assume $\bx = 1$ without
loss of generality. Let $\pz_0
\sm \pz_1 \sm \dots \sm \pz_\NNN$ be the increasing
enumeration of¨$\Div(\pz)$. By definition, all prefixes
of¨$\ww$ represent divisors of¨$\pz$, so, letting¨$\ell$
be the length of¨$\ww$, there exists a mapping $f : \{0,
\dots, \ell\} \to  \{0,
\dots, \NNN\}$ such that, for each¨$\kk$, the
length¨$\kk$ prefix of¨$\ww$
represents¨$\pz_{f(\kk)}$. By construction, we have
$f(0) = 0$ and $f(\ell) = \NNN$.

The difference with Lemma~\ref{L:Enum} is
that $f$ need not be increasing. Now, let
$p_1$, \dots, $p_{\NN_\mm}$ be the
$\NN_\mm$~positions in¨$\ww$ where the
generator~$\ss m$ occurs, completed with
$p_0 = 0$. Then, in the prefix of¨$\ww$ of
length¨$p_1$, \ie, in the subword of¨$\ww$
corresponding to positions from¨$p_0+1$
to¨$p_1$, there is one¨$\ss m$, plus
letters¨$\ss i^{\pm1}$ with $i <
m$ (Figure¨\ref{F:Jump}). This subword is
therefore $\s$-positive, hence we
must have $\pz_{f(p_0)} \sm \pz_{f(p_1)}$,
which requires
$f(p_0) < f(p_1)$. Moreover, the quotient
$\pz_{f(p_0)}\inv
\pz_{f(p_1)}$ is a braid that admits at least
one $\s$-positive expression containing¨$\ss m$,
hence $\pz_{f(p_0)}
\sm_\mm \pz_{f(p_1)}$ holds. Now the same
is true between¨$f(p_1)$ and¨$f(p_2)$,
\etc\ Hence the number of $\mm$-jumps in the
increasing enumeration of~$\Div(\pz)$ is at
least~$\NN_\mm$, \ie, we have $\hh\mm(\pz) \ge
\NN_\mm + 1$.
\end{proof}

It remains to prove the second inequality in
Proposition~\ref{P:MainHeight}, \ie, to prove
that, if $\pz$ is a positive $\nn$-braid
satisfying $\hh{\mm}(\pz) = \NN+1$, then
$\pz$ admits an almost-positive expression
containing $\NN$~generators~$\ss\mm$. 
The problem is as follows: if $\pz$ is a
positive braid and $\bx, \by$ are left
divisors of~$\pz$ satisfying $\bx \sm \by$,
then, by definition, the quotient $\bx\inv
\by$ admits some $\s$-positive
expression~$\ww$, but nothing {\it
a priori } guarantees that $\ww$ be drawn
in~$\CG(\pz)$. In other words, we might have
$\bx \sm \by$ but no $\s$-positive witness
for this inequality inside~$\Div(\pz)$.
This however cannot happen, but the proof
requires a rather delicate argument. 

\begin{prop} \label{P:Handle}
Let $\pz$ be a positive braid. Then, 
for all $\bx, \by$ in¨$\Div(\pz)$,
the following are equivalent:

$(i)$ The relation $\bx \sm \by$ holds,
\ie, there exists a $\s$-positive path
from~$\bx$ to~$\by$ in the Cayley graph
of~$B_\infty$;

$(ii)$ There exists a $\s$-positive path 
from¨$\bx$ to¨$\by$ in the Cayley graph
of~$B_\nn$;

$(iii)$ There exists a $\s$-positive path 
from¨$\bx$ to¨$\by$ in¨$\CG(\pz)$.
\end{prop}

\begin{proof}
Clearly $(iii)$ implies¨$(ii)$, which 
implies¨$(i)$. We shall prove that $(i)$
implies~$(iii)$---and thus reprove that
$(i)$ implies~$(ii)$, which was first
proved in¨\cite{Lrb}---by
using the handle reduction method of¨\cite{Dfo, Dgr}. The
problem is to prove that, among all $\s$-positive paths
connecting¨$\bx$ to¨$\by$ in the Cayley graph
of¨$B_\infty$, at least one is drawn in~$\CG(\pz)$.

Now, let $\uu, \vv$ be positive words representing
$\bx$ and¨$\by$. Then the word¨$\uu\inv\vv$
represents¨$\bx\inv\by$, and, by hypothesis, it is drawn
 in~$\CG(\pz)$ from¨$\bx$. Handle
reduction is an operation that transforms a braid word into
equivalent words and eventually produces a $\s$-positive
word if it exists. It is proved in¨\cite{Dfo} that, for
every $\nn$¨strand braid word¨$\ww$, there exists a
finite fragment¨$\G_\ww$ of the Cayley
graph of~$B_\nn^+$ and a vertex¨$\bx_\ww$
of¨$\G_\ww$ such that
$\ww$ and all words obtained from¨$\ww$ by handle reduction
are drawn from¨$\bx_\ww$ in¨$\G_\ww$. Moreover,
when $\ww$ has the form¨$\uu\inv\vv$ with $\uu, \vv$
positive, then all vertices in¨$\G_\ww$ are the
left divisors of the least common right multiple of
the braids represented by¨$\uu$ and¨$\vv$, here
¨$\bx$ and¨$\by$,  while $\bx_\ww$ is the braid
represented by¨$\uu$, \ie,¨$\bx$. As
$\bx$ and¨$\by$ are divisors of¨$\pz$, so is their
least common right multiple, and the
graph¨$\G_\ww$ is included in~¨$\CG(\pz)$.
It follows that every word obtained
from¨$\uu\inv\vv$ using handle reduction is
drawn from¨$\bx$ in¨$\CG(\pz)$. The
termination of handle reduction guarantees
that, among these words, at least one is
$\s$-positive, so $(iii)$ follows.
\end{proof}

A direct application of Proposition~\ref{P:Handle} is the
existence of $\s$-positive quotient-sequences drawn in
the Cayley graph. The definition is as follows:

\begin{defi} \label{D:Witness}
Assume that $\pz$ is a positive braid and
$\XX$ is a subset of~$\Div(\pz)$. Let
$\bx_0 \sm \dots \sm \bx_\NN$ be the
increasing enumeration of~$\XX$. We say that
a sequence of words $\vec\ww = (\ww_1, \dots,
\ww_\NNN)$ is a {\it quotient-sequence} for~$\XX$ if,
for each~$\kk$, the word~$\ww_\kk$ is an  expression
of~$\bx_{\kk-1}\inv \bx_\kk$ for each¨$\kk$. We
say that $\vec\ww$ is {\it $\s$-positive} if every
entry in¨$\vec\ww$ is $\s$-positive, and that
$\vec\ww$ is {\it drawn in¨$\CG(\pz)$} (from¨$\bx_0$)
if $\ww_\kk$ is drawn from~$\bx_{\kk-1}$ in~$\CG(\pz)$
for each¨$\kk$.
\end{defi}

\begin{coro} \label{C:Witness}
Assume that $\pz$ is a positive braid. Then
every subset of~$\Div(\pz)$ admits a
$\s$-positive quotient-sequence drawn in¨$\CG(\pz)$.
\end{coro}

\begin{figure} [htb]
\begin{picture}(60,45)(0, -0)
\put(0,0){\includegraphics{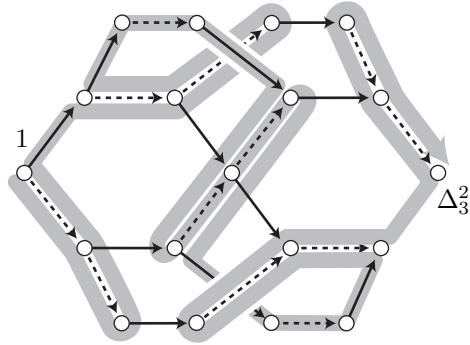}}
\put(1,26){$1$}
\put(57,18){$\D_3^2$}
\end{picture}
\caption{\smaller The increasing enumeration of the
divisors of¨$\DD32$, and a $\s$-positive
quotient-sequence drawn in¨$\GG32$: the associated path
visits every vertex, and is labeled
$\tt aaAAbaaAAbabAAaaAAbabAAaabAAaa$; it crosses
$6$ $\ss2$-edges (and no $\ss2\inv$)}
\label{F:MaxPath}
\end{figure}

\begin{exam} \label{X:Jump}
 (Figure¨\ref{F:MaxPath})
By computing the successive quotients in the
increasing enumeration of¨$\DDD32$ given in
Example¨\ref{X:List}, we easily find that
$$\tt (a, a, AAb, a, a, AAb, a, AAb, a, a, bAA, a,
bAA, a, a, bAA, a, a)$$
is a $\s$-positive quotient-sequence for $\Div(\DD32)$
drawn in¨$\GG32$. This sequence turns out to be the
unique sequence with the above properties, but this
uniqueness is specific to the case of $3$-braids ({\it
cf.} Figure~\ref{F:Decomp41} below). 
\end{exam}

We can now easily complete the proof of
Proposition~\ref{P:MainHeight}:

\begin{proof}[Proof of
Proposition~\ref{P:MainHeight}]
Let $(\pz_0, \dots, \pz_\NN)$
be the $\sm$-increasing enumeration\break
of $\Div(\pz)$. By Corollary¨\ref{C:Witness},
there exists a $\s$-positive quotient-sequence¨$\vec\ww$
for¨$\Div(\pz)$ that is drawn in¨$\CG(\pz)$. Let $\ww =
\ww_1 \ldots \ww_\NN$. By construction, $\ww$ is a
$\s$-positive word drawn in¨$\CG(\pz)$, and the number of
occurrences of the main generator¨$\ss\mm$ in¨$\ww$ is
(at least) the number of $\mm$-jumps in¨$(\pz_0, \dots,
\pz_\NN)$. So we have $\cc(\pz) \ge \hh\mm(\pz)
-1$. Owing to Lemma¨\ref{L:MainHeight}, this
completes the proof.
\end{proof}

\begin{rema}
Assume that $\vec\ww$ is a $\s$-positive
quotient-sequence  for~$\Div(\pz)$, and
$\ss\mm$ is the main generator
occurring in~$\vec\ww$. Then
each word~$\ww_i$ contains zero
or one letter~$\ss\mm$. Indeed, if $\ww_i$
contained two~$\ss\mm$'s or more, then the
vertex reached after the first~$\ss\mm$ ought
to lie in the open $\sm$-interval determined by two
successive entries of~$\vec\pz$, and the
latter is empty by construction since all
elements of~$\Div(\pz)$ occur in~$\vec\pz$.
\end{rema}

\section{A decomposition result for~$\DDiv\pz$}
\label{S:Deter}

In this section, we establish a structural result
describing $(\DDD\nn\dd, <)$ as the concatenation of
$\ccc\nn\dd+1$ intervals isomorphic to subsets
of¨$(\DDD{\nn-1}\dd, <)$. We deduce an explicit formula
connecting¨$\hhh\rr\nn\dd$ with the number of braids
in¨$\DDD\nn\dd$ whose $\dd$th factor is right divisible
by¨$\D_\rr$, which in turn enables us to complete the
computation of¨$\ccc\nn\dd$ and¨$\hhh\rr\nn\dd$ for
small values of¨$\rr$, $\nn$ and¨$\dd$.

\subsection{$B_\rr$-classes}

In order to analyse the linearly ordered
sets¨$\DDDD\nn\dd$, and, more generally, 
$(\Div(\pz), \sm\nobreak)$ for¨$\pz$ a positive braid,
we introduce convenient partitions. As $B_\rr$ is a
group for each~$\rr$, it is clear that the relation
$\bx\inv \by \in B_\rr$ defines an equivalence
relation on (positive) braids, so we may put:

\begin{defi}
For $\rr \ge 1$ and $\bx, \by$
in¨$B_\infty^+$, we say that $\bx$ and
$\by$ are {\it $B_\rr$-equivalent} if
$\bx\inv \by$ belongs to¨$B_\rr$.
\end{defi}

By construction, $B_\rr$-equivalence is
compatible with multiplication on the left.
In the sequel, we consider the restriction
of $B_\rr$-equivalence to finite subsets
of~$B_\infty^+$ of the form~$\Div(\pz)$,
\ie, we use $B_\rr$-equivalence to partition
$\Div(\pz)$ into subsets, naturally called
$B_\rr$-classes.

\begin{figure} [htb]
\begin{picture}(100,43)(0, -0)
\put(0,0){\includegraphics{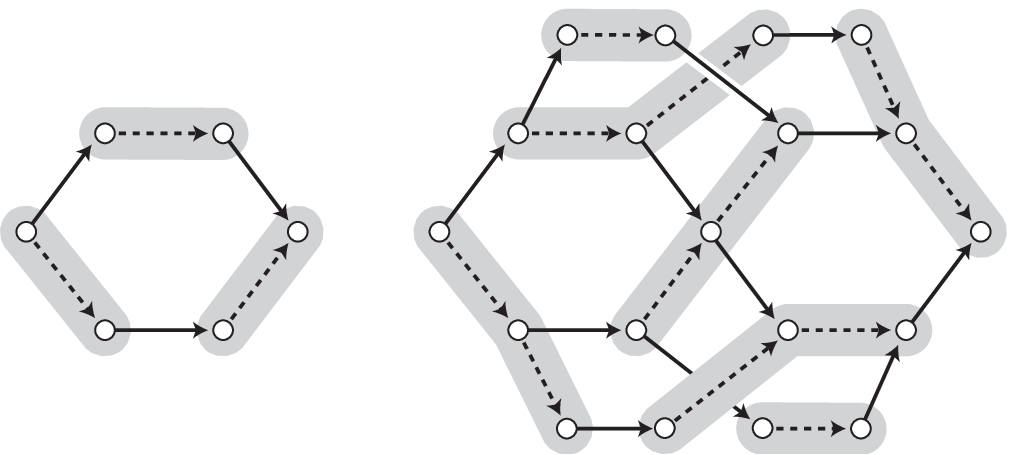}}
\end{picture}
\caption{\smaller The $B_2$-classes
in¨$\DDD3{}$ and¨$\DDD32$}
\label{F:Classes}
\end{figure}

\begin{exam}
As $B_1$ is trivial, $B_1$-equivalence is
equality, and, therefore, the $B_1$-classes are
singletons. On the other hand, any two elements
of¨$B_\nn$ are
$B_\rr$-equivalent for each¨$\rr \ge n$, so, for $\pz$
in¨$B_\nn^+$, there is only one $B_\rr$-class for $\rr \ge
n$, and the only interesting relations arise for
$1 < \rr < n$. For instance, $\DDD3{}$
contains three $B_2$-classes, while $\DDD32$
contains seven of them
(Figure¨\ref{F:Classes}).
\end{exam}

Saying that there is an $\rr$-jump
between two braids¨$\bx$ and¨$\by$ means that
$\bx\inv \by$ is $\s$-positive and does not
belong to¨$B_\rr$, so, for $\bx \sm \by$, we
have the equivalence 
\begin{equation} \label{E:Jump}
\big(\mbox{$\bx, \by$ are not
$B_\rr$-equivalent}\big)
\Longleftrightarrow
\left(
\begin{matrix}
\mbox{there is a $\rr$-jump between} \\
\mbox{between¨$\bx$ and¨$\by$}
\end{matrix}
\right).
\end{equation} 

\begin{lemm} \label{L:Interval}
Assume that $\pz$ is a positive braid. Then,
each $B_\rr$-class in¨$\Div(\pz)$ is an
interval for¨$\sm$ and there is an $\rr$-jump
between each $B_\rr$-class and the next one.
\end{lemm}

\begin{proof}
Assume $\bx \sm \by \in \Div(\pz)$. By
\eqref{E:Jump}, if $\bx$ and $\by$
are not $B_\rr$-equivalent, there is an
$\rr$-jump between¨$\bx$ and¨$\by$, hence
between¨$\bx$ and any element of¨$\Div(\pz)$
above¨$\by$, so no such element may be
$B_\rr$-equivalent to¨$\bx$. This implies
that each $B_\rr$-class is an
$\sm$-interval.
\end{proof}

\begin{coro} \label{C:NbClasses}
For each $\rr \ge 1$, 
the number of $B_\rr$-classes in¨$\Div(\pz)$
is¨$\hh\rr(\pz)$.
\end{coro}

\begin{proof}
By¨\eqref{E:Jump}, there is no
$\rr$-jump between two elements of the same
$B_\rr$-class, and there is one between two
elements not in the same $B_\rr$-class. Thus
the number of $B_\rr$-classes is the number
of $\rr$-jumps in the $\sm$-increasing
enumeration of¨$\Div(\pz)$ augmented by¨$1$,
hence, by Lemma~\ref{L:Enum}, it is¨$\hh\rr(\pz)$.
\end{proof}

$B_\rr$-equivalence provides a partition
of¨$(\Div(\pz), \sm)$ into finitely many
subintervals. The interest of this partition
is that we can describe $B_\rr$-classes rather
precisely and, typically, connect them with
subsets of¨$B_\rr$. In particular, this will
allow for connecting the ordered
sets¨$\DDDD\nn\dd$ with the
sets¨$(\DDD{\nn-1}\dd, \sm)$.

\begin{figure} [htb]
\begin{picture}(86,44)(0, 0)
\put(2,0){\includegraphics{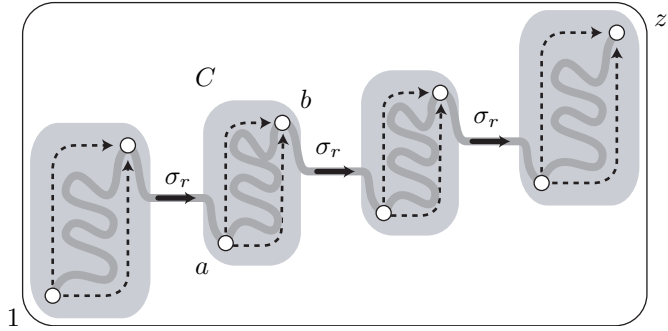}}
\put(25,32){$C$}
\put(25,7){$a$}
\put(39,29){$b$}
\put(0,0){$1$}
\put(21,19){$\ss\rr$}
\put(41,23){$\ss\rr$}
\put(62,27){$\ss\rr$}
\put(86,40){$\pz$}
\end{picture}
\caption{\smaller Decomposition of $(\Div(\pz),
\sm)$ into $B_\rr$-classes: each class¨$C$ is a
lattice with respect to divisibility; the increasing
enumeration of¨$\Div(\pz)$ exhausts the first class,
then jumps to the next one by an $\rr$-jump, \etc\ The
number of classes is¨$\hh\rr(\pz)$.}
\label{F:Decomp}
\end{figure}

\begin{prop} \label{P:Lattice}
(Figure¨\ref{F:Decomp})
Assume $\pz \in B_\infty^+$ and $\rr \ge 1$. Let
$\Class$ be a
$B_\rr$-class in¨$\Div(\pz)$, and let
$\ea, \eb$ be its $\sm$-extremal elements. Then $\ea$
divides every element of¨$\Class$ on the left, and the left
translation by¨$\ea$ defines an isomorphism of
$(\Div(\ea\inv\eb), \dive, \sm)$ onto¨$(\Class, \dive,
\sm)$. In particular, $(\Class, \dive)$ is a
lattice.
\end{prop}

\begin{proof}
By Lemma¨\ref{L:Interval}, $\Class$ is the
$\sm$-interval determined by¨$\ea$ and¨$\eb$,
\ie, we have $$\Class = \{\bx
\in \Div(\pz) ; \ea \sm \bx \sm \eb\}.$$

We know that $\Div(\pz)$ is a
lattice with respect to left divisibility:
any two elements¨$\bx, \by$ of¨$\Div(\pz)$
admit a greatest left common divisor, here
denoted $\gcd(\bx, \by)$, and a least common
right multiple, denoted¨$\lcm(\bx, \by)$.
Firstly, we claim that $\Class$ is a lattice with
respect to left divisibility, \ie, the left gcd 
and the right lcm of two elements of¨$\Class$
lie in¨$\Class$. So assume $\bx, \by \in \Class$. Let $\bx_0,
\by_0$ be defined by $\bx = \gcd(\bx, \by)
\bx_0$ and $\by = \gcd(\bx,
\by) \by_0$.  The hypothesis that $\bx\inv\by$ belongs
to¨$B_\rr$ implies that there exist $\bx_1, \by_1$ in¨$B_\rr^+$
satisfying $\bx\inv \by = \bx_1\inv \by_1$. By definition of
the gcd, there must exist a positive
braid¨$\pz_1$ satisfying $\bx_1 = \pz_1 \bx_0$
and $\by_1 = \pz_1 \by_0$. Because
$\pz_1$ is positive, $\bx_1 \in B_\rr^+$
implies $\bx_0 \in B_\rr^+$, hence
$\gcd(\bx, \by) \in \Class$. As for the lcm,
the conjunction of $\bx = \gcd(\bx,
\by)
\bx_0$ and $\by = \gcd(\bx, \by) \by_0$
implies $$\lcm(\bx, \by) = \gcd(\bx,
\by)\lcm(\bx_0, \by_0).$$ As $\bx_0, \by_0
\in B_\rr^+$ implies $\lcm(\bx_0, \by_0) \in
B_\rr^+$, we deduce $\lcm(\bx, \by) \in \Class$.

As $\Class$ is finite, it follows that $\Class$
admits a global gcd. Because the linear
ordering¨$\sme$ extends the partial
divisibility ordering¨$\dive$, this global
gcd must be the $\sm$-minimum¨$\ea$ of¨$\Class$.
Symmetrically, $\Class$ admits a global lcm,
which must be the $\sm$-maximum¨$\eb$. So, at
this point, we know that $\ea$ is a left
divisor of every element in¨$\Class$, and $\eb$
is a right multiple of each such element,
\ie, we have 
\begin{equation} \label{E:Inclusion}
\Class \ince \{\bx \in B_\infty^+ ; \ea \dive
\bx \dive \eb\}.
\end{equation}
Moreover, $\ea \dive \bx \dive \eb$
implies  $\ea \sme \bx \sme \eb$, hence $\bx
\in \Class$, so the inclusion
in¨\eqref{E:Inclusion} is an
equality.

Now, put $F(\bx) = \ea\bx$ for¨$\bx$
in¨$\Div(\ea\inv \eb)$. As $B_\infty^+$ is
left cancellative, $F$ is injective.
Moreover, for $\bx$ a positive braid, $\bx
\dive
\ea\inv
\eb$ is equivalent to $\ea\bx \dive \eb$, so
the image of¨$F$ is $\{\bx \in B_\infty^+ ;
\ea \dive
\bx \dive \eb\}$, hence is¨$\Class$. Finally, by
construction, $F$ preserves both¨$\dive$
and¨$\sm$.
\end{proof}

For¨$\rr = 1$, each $B_\rr$-class is a singleton, and
Proposition¨\ref{P:Lattice} says nothing; similarly, if
the main generator of¨$\pz$ is¨$\ss\mm$,
there is only one $B_\rr$-class for $\rr >
\mm$, and we gain no information. But, for
$1 < \rr \le \mm$, and specially for $\rr =
\mm$, Proposition¨\ref{P:Lattice} states that
the chain¨$\Div(\pz)$ is obtained by
concatenating $\hh\rr(\pz)$ copies of sets of
the form¨$\Div(\pz')$ with¨$\pz'$ of index at
most¨$\rr$. In particular, for $\pz = \DD
n\dd$, we have

\begin{coro} \label{C:Structure}
For each¨$\nn$ and each¨$\rr$ with $\rr < \nn$, the
chain $\DDDD\nn\dd$ is obtained by  concatenating
$\hhh\rr\nn\dd$ in\-tervals, each of which,
when equipped with¨$\dive$, is a translated
copy of some initial sublattice
of¨$(\DDD\rr\dd, \dive)$.
\end{coro}

\begin{figure} [htb]
\begin{picture}(64,43)(0, 0)
\put(2,0){\includegraphics{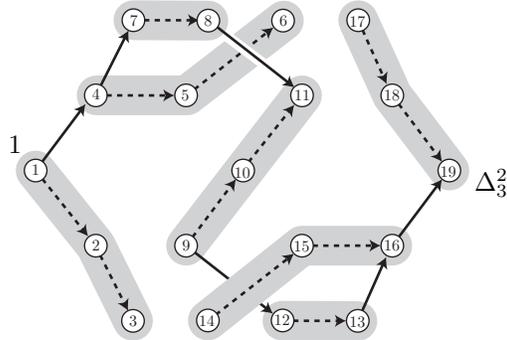}}
\put(1,25){$1$}
\put(63,20){$\DD32$}
\end{picture}
\caption{\smaller Decomposition of
$\DDDD32$ into $B_2$-classes: the increasing
enumeration of¨$\DDDD32$ is the concatenation
of the increasing enumeration of the successive
classes, separated by $2$-jumps (compare with
Figure¨\ref{F:MaxPath}); in this case, $B_2$-classes
are simply chains with respect to divisibility}
\label{F:Decomp32}
\end{figure}

\begin{figure} [htb]
\begin{picture}(104,48)(0, 0)
\put(2,0){\includegraphics{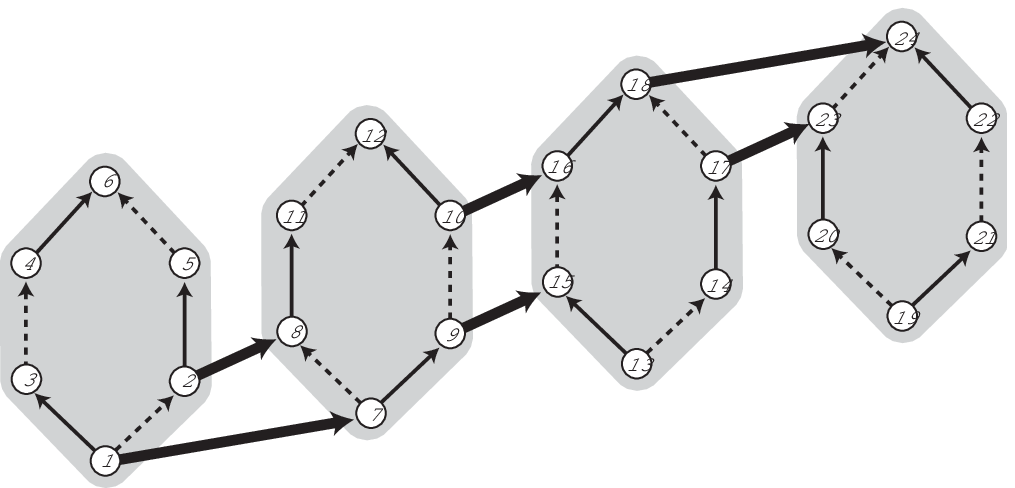}}
\put(7,0){$1$}
\put(97,46){$\DD4{}$}
\end{picture}
\caption{\smaller Decomposition of
$\DDDD4{}$ into $B_3$-classes; note that the
$\ss3$-arrows (thick) corresponding to 3-jumps
are not unique; in this case, all $B_3$-classes are
isomorphic to the lattice $(\DDD3{}, \sm, \dive)$, \ie,
to the Cayley graph of¨$\D_3$}
\label{F:Decomp41}
\end{figure}

The case of¨$\DD32$ and¨$\DD4{}$ are illustrated in
Figures¨\ref{F:Decomp32} and¨\ref{F:Decomp41}.

\subsection{Extremal elements}

The next step is to observe that extremal
points in $B_\rr$-classes admit a simple
characterization in terms of divisibility.

\begin{prop} \label{P:Extremal}
Assume that $\pz$ is a positive braid. 

$(i)$ An element¨$\bx$ of¨$\Div(\pz)$ is the maximum of its
$B_\rr$-class if and only if the relation $\bx\ss i
\dive\pz$ fails for $1 \le i < \rr$.

$(ii)$ An element¨$\bx$ of¨$\Div(\pz)$ is the minimum of its
$B_\rr$-class if and only if no¨$\ss i$ with $1 \le i <
\rr$ divides¨$\bx$ on the right.
\end{prop}

\begin{proof}
$(i)$ The condition is necessary: if $\bx \ss i$ lies
in¨$\Div(\pz)$ for some¨$i$ with $i < \rr$, then $\bx \ss
i$ lies in the same $B_\rr$-class as¨$\bx$, and it is larger
both for¨$\dive$ and¨$\sm$, so $\bx$ cannot be maximal in
its $B_\rr$-class. Conversely, assume that $\bx$ is
not maximal in its $B_\rr$-class. Then there
exists¨$\by$ satisfying $\bx \sm \by$ and $\by$ is
$B_\rr$-equivalent to¨$\bx$. Now, by 
Proposition¨\ref{P:Lattice},
the lcm of¨$\bx$ and¨$\by$ is also $B_\rr$-equivalent
to¨$\bx$, which means that there exists¨$\by_1$ in¨$B_\rr^+$
satisfying $\lcm(\bx, \by) = \bx\by_1$. Now $\bx \sm \by$
implies $\by_1 \not= 1$, so there must exist¨$i < m$ such
that $\ss i$ is a left divisor of¨$\by_1$. Then we
have $\bx\ss i \dive \bx\by_1 \dive \pz$, hence $\bx
\ss i \dive \pz$.

$(ii)$ The argument is symmetric. If we have $\bx = 
\by \ss i$ for some positive braid¨$\by$ and $i < \rr$, then
$\by$ belongs to the $B_\rr$-class of¨$\bx$, and $\bx$
cannot be minimal in its $B_\rr$-class. Conversely, assume
that $\bx$ is not minimal in its $B_\rr$-class. Then there
exists¨$\by$ satisfying $\by \sm \bx$ and $\by$ is
$B_\rr$-equivalent to¨$\bx$. By
Proposition¨\ref{P:Lattice} again, the gcd
of¨$\bx$ and¨$\by$ is also $B_\rr$-equivalent
to¨$\bx$, which means that there
exists¨$\by_0$ in¨$B_\rr^+$ satisfying
$\gcd(\bx, \by) \by_0 = \bx$. As $\by \sm
\bx$ implies $\by_0 \not= 1$, there must exist¨$i < m$
such that $\ss i$ is a right divisor of¨$\by_0$, hence
of¨$\bx$.
\end{proof}

When we apply the previous criterion to the
braids¨$\DD\nn\dd$, we obtain:

\begin{prop} \label{P:MaxClass}
For $\bx$ in¨$\DDD\nn\dd$ and $1 \le
\rr \le \nn$, the
following are equivalent:

$(i)$ The element¨$\bx$ is $\sm$-maximal in its
$B_\rr$-class;

$(ii)$ The element¨$\bx \ss i$ belongs to¨$\Div(\DD\nn\dd)$
for no¨$i < \rr$;

$(iii)$ The $\dd$th factor of¨$\bx$ is
right divisible by¨$\D_\rr$.

$(iv)$ The $\dd+1$st factor of¨$\bx \D_\rr$
is¨$\D_\rr$.
\end{prop}

\begin{proof}
The equivalence of¨$(i)$ and¨$(ii)$ is given
by Proposition¨\ref{P:Extremal}$(i)$. It
remains to establish the equivalence
of¨$(ii)$--$(iv)$. For $\rr = 1$, $(ii)$ is
vacuously true, while $(iii)$ and
$(iv)$ always hold. So the expected equivalences are true.
We henceforth assume $\rr \ge 2$.

Let $\bx$ belong to¨$\DDD\nn\dd$, and let $\sx_\dd$ be the $\dd$th
factor in the normal form of¨$\bx$. For $i < n$, saying that
$\bx \ss i$ does not belong to¨$\Div(\DD\nn\dd)$ means that
the normal form of¨$\bx \ss i$ has length $\dd+1$, hence,
equivalently, that the normal form of¨$\sx_\dd \ss i$ has
length¨$2$. This occurs if and only if
$\ss i$ is a right divisor of¨$\sx_\dd$. So, for
$\rr \le n$,
$(ii)$ is equivalent to $\sx_\dd$ being right divisible by
all¨$\ss i$'s with $1 \le i < \rr$, hence to
$\sx_\dd$ being right divisible by the (left)
lcm of these elements, which is¨$\D_\rr$.

Finally, $(iii)$ and $(iv)$ are equivalent. Indeed, if
the $\dd$th factor¨$\sx_\dd$ in the normal form of¨$\bx$ is
divisible by¨$\D_\rr$ on the right, then $(\sx_\dd, \D_\rr)$ is a
normal sequence as no $\ss i$ with $i < \rr$
from¨$\D_\rr$ may pass to¨$\sx_\dd$. Hence
$(\sx_1, \ldots, \sx_\dd, \D_\rr)$ is a normal
sequence, necessarily the normal form
of¨$\bx\D_\rr$.  Conversely, assume that the
normal form of¨$\bx\D_\rr$ is
$(\sx_1, \ldots, \sx_\dd, \D_\rr)$. The hypothesis that
$(\sx_\dd, \D_\rr)$ is normal implies that $\sx_\dd$ is divisible
on the right by each¨$\ss i$ with $i < \rr$,
hence is divisible on the right by¨$\D_\rr$.
Now $(\sx_1, \ldots,
\sx_\dd)$ is the normal form of¨$\bx$.
\end{proof}

Observe that, for $\rr \ge 2$, an element
of¨$\DDD\nn\dd$ that is $\sm$-maximal in its
$B_\rr$-class cannot
belong to¨$\DDD\nn{\dd-1}$, \ie, cannot have
degree $\dd-1$ or less, since the
$\dd$th factor of its normal form cannot be¨$1$.

Similar conditions characterize the minimal elements of the
$B_\rr$-classes. Because the normal form has a priviledged
orientation, the results are not entirely symmetric of
those of Proposition¨\ref{P:MaxClass}

\begin{prop} \label{P:MinClass}
For $\bx$ in¨$\DDD\nn\dd$ and $1 \le
\rr \le \nn$, the following are equivalent:

$(i)$ The element¨$\bx$ is $\sm$-minimal of its
$B_\rr$-class;

$(ii)$ No¨$\ss i$ with $i < \rr$ is a right divisor
of¨$\bx$;

$(iii)$ The degrees of¨$\bx$ and¨$\bx \D_\rr$ are equal.
\end{prop}

\begin{proof}
The equivalence of¨$(i)$ and¨$(ii)$ is given by
Proposition¨\ref{P:Extremal}$(ii)$. On the other hand,
everything is obvious for¨$\rr = 1$. So it remains to
establish the equivalence of¨$(ii)$ and¨$(iii)$ in the case
$\rr \ge 2$. Now, assume that¨$(ii)$ holds an $\bx$ has
degree¨$\dd$. The hypothesis that $\ss i$ is not a right
divisor of¨$\bx$ implies that $\bx \ss i$ is a divisor
of¨$\DD\nn\dd$. As this holds for each $i <
\rr$, the lcm of $\bx \ss1$, \dots, $\bx \ss{\rr-1}$, which
is $\bx \D_\rr$, also divides¨$\DD\nn\dd$, which means that
$\bx \D_\rr$ has degree (at most)¨$\dd$. So $(ii)$
implies¨$¨(iii)$.

Conversely, assume that $\ss i$ divides¨$\bx$ on the right.
Then the degree of¨$\bx \ss i$ is strictly larger than that
of¨$\bx$, and, {\it a fortiori}, the same is true for¨$\bx
\D_\rr$.
\end{proof}

\subsection{Determination of¨$\hhh\rr\nn\dd$}

A direct application of the previous results is a formula
connecting the number of $B_\rr$-classes
in¨$\Div(\DD\nn\dd)$, \ie, the
numbers¨$\hhh\rr\nn\dd$, with the number of
braids whose normal form ends with some
specific factor.

\begin{defi}
For $\nn, \dd \ge 1$ and for¨$\sp$ a simple
$\nn$-braid, we denote by¨$\bbb\nn\dd\sp$ the number
of positive braids of degree at most¨$\dd$, \ie, of
divisors of¨$\DD\nn\dd$, whose $\dd$th factor is¨$\sp$.
\end{defi}

\begin{prop} \label{P:Main}
For $1 \le \rr \le \nn$, we have
\begin{equation} \label{E:Main}
\hhh\rr\nn\dd = \sum_{\mbox{\Small $s$ right divisible by
$\D_\rr$}} \bbb\nn\dd s  = \bbb \nn{\dd+1}{\D_\rr} .
\end{equation} 
\end{prop}

In words: The number of $\rr$-jumps in $\DDDD\nn\dd$ is
the number of
$\nn$-braids of degree at most¨$\dd$ whose $\dd$th factor
is right divisible by¨$\D_\rr$.

\begin{proof}
By Corollary¨\ref{C:NbClasses},
$\hhh\rr\nn\dd$ is the number of
$B_\rr$-classes in¨$\Div(\DD\nn\dd)$. Each
class contains exactly one maximum element, and, by
Proposition¨\ref{P:MaxClass}, the latter are
characterized by the property that their
$\dd$th factor is right
divisible by¨$\D_\rr$. The first equality
in¨\eqref{E:Main} follows. The second one
follows from the equivalence of¨$(iii)$
and¨$(iv)$ in Proposition¨\ref{P:MaxClass}.
\end{proof}

For $\rr = 1$, as every simple braid is divisible by¨$1$ on
the right, Relation¨\eqref{E:Main} reduces to
\begin{equation} \label{E:MainOne}
\hhh1\nn\dd = \sum_{\sp} \bbb\nn\dd\sp = \bbb
\nn{\dd+1}{1},
\end{equation}
a special case of the relation $\hh1(\pz) =
\card\Div(\pz)$ of
Proposition¨\ref{P:Decreasing}. 
For $\rr = \nn$, as the only normal sequence of
length¨$\dd$  that finishes with¨$\D_\nn$ is  $(\D_\nn,
\dots, \D_\nn)$, Relation¨\eqref{E:Main}
reduces to
\begin{equation} \label{E:Mainn}
\hhh\nn\nn\dd = 1,
\end{equation}
also noted in Proposition¨\ref{P:Decreasing}.
Finally, for $\rr = \nn-1$, we obtain using
Proposition¨\ref{P:MainHeight}:

\begin{coro} \label{C:Main}
For $\nn \ge 2$, we have
\begin{equation} \label{E:MainNO}
\ccc\nn\dd = \hhh{\nn - 1}\nn\dd -1 = \sum_{i =
2}^{\nn}
\bbb\nn\dd{\ss i \ss{i+1} \dots \ss{\nn-1}\D_{\nn - 1}} =
\bbb\nn{\dd+1}{\D_{\nn - 1}} - 1.
\end{equation}
\end{coro}

\begin{proof}
The simple $\nn$-braids that are right divisible
by¨$\D_{\nn-1}$ are the braids of the form
$\ss i \ss{i+1} \dots \ss{\nn-1}$ with $1 \le i\le
\nn$. Indeed,
it is clear that every such braid is simple
and right divisible by¨$\D_{\nn-1}$.
Conversely, the only possibility for $\pz
\D_{\nn-1}$ to be simple is that $\pz$ moves the $\nn$th
strand to some position between¨$1$ and¨$\nn$, but
introduces no crossing between the remaining strands.
Finally, $\ss1 \ss2 \dots \ss{\nn-1}\D_{\nn-1}$
is¨$\D_\nn$, and we already observed that
$\bbb\nn\dd{\D_\nn}$ is¨$1$, so we obtain the first
equality in¨\eqref{E:MainNO}.
\end{proof}

\subsection{Computation of¨$\bbb\nn\dd\sp$}

By Lemma¨\ref{L:Normal}, normal sequences are
characterized by a local condition involving only pairs
of consecutive elements. It follows that the set of all
normal sequences is a rational set, \ie, it can be
recognized by a finite state automaton. Standard
arguments then show that the numbers¨$\bbb\nn\dd\sp$ obey
a linear recurrence. Building on this observation,
seemingly first used in the case of braids in¨\cite{Chb},
we can obtain explicit formulas for the
parameters¨$\ccc\nn\dd$ and¨$\hhh\rr\nn\dd$ for small
values of¨$\rr$, $\nn$, and/or¨$\dd$. We shall not go
into details in the current paper, and refer
to¨\cite{Dhi} where all formulas are established---and
where, more generally, the rich combinatorics underlying
the normal form of braids is investigated.

In the sequel, we write $(M)_{\xx, \yy}$ for the $(\xx,
\yy)$-entry of a matrix¨$M$. The general principle for
computing the numbers¨$\bbb\nn\dd\sp$ for some
fixed¨$\nn$ is as follows:

\begin{lemm} \label{L:Comput}
For $\nn \ge 1$, let $M_\nn$ be the square matrix with
entries indexed by simple $\nn$-braids defined by
$$
(M_\nn)_{\sp, \spp} = 
\begin{cases}
1 & \mbox{if $(\sp, \spp)$ is normal,}\\
0 & \mbox{otherwise.}
\end{cases}$$
Then, for every simple¨$\spp$ and every
$\dd \ge 1$, we have $\bbb\nn\dd{\spp} = ((1, 1, \dots, 1)
\, M_\nn\!\!\!{}^{\dd-1})_\spp$.
\end{lemm}

The proof is an easy induction on¨$\dd$.

\begin{exam}
The matrix¨$M_1$ is
$(1)$, corresponding to
$\bbb1\dd1 = 1$. For $\nn=2$, using the enumeration $(1,
\ss1)$ of simple $2$-braids, we find $M_2 = 
\left(\begin{matrix}1&0\\1&1\end{matrix}\right)$,
leading to
$\bbb2\dd1 = \dd$, $\bbb2\dd{\ss1} = 1$, as could
be expected: there are $\dd+1$ braids of degree
at most¨$\dd$, namely the braids¨$\ss1^\expo$ with
$\expo < \dd$, whose $\dd$th factor is¨$1$, and
$\ss1^\dd$, whose $\dd$th factor is¨$\D_2$, \ie,
¨$\ss1$. For $\nn = 3$, using the enumeration $(1, \ss1,
\ss2,\ss2\ss1, \ss1\ss2, \D_3)$ of simple $3$-braids,  we
obtain
$$M_3 = \left(
\begin{matrix}
1&0&0&0&0&0\\
1&1&0&0&1&0\\
1&0&1&1&0&0\\
1&1&0&0&1&0\\
1&0&1&1&0&0\\
1&1&1&1&1&1
\end{matrix}
\right),$$
from which we can deduce for instance $\bbb331= 19$ or
$\bbb34{\ss1} = 15$ using Lemma¨\ref{L:Comput}.
\end{exam}

Using Proposition¨\ref{P:Main}, we deduce

\begin{prop} \label{P:Values}
With $M_\nn$ as in Lemma¨\ref{L:Comput}, we have
for $\nn \ge \rr \ge 1$ and $\dd \ge 1$
\begin{gather*}
\ccc\nn\dd = ((1, 1, \dots, 1) \, M_\nn^\dd)_{\D_{\nn-1}}
- 1,\\
\hhh\rr\nn\dd = ((1, 1, \dots, 1) \,
M_\nn^\dd)_{\D_{\rr}}.
\end{gather*}
\end{prop}

\begin{coro}
$(i)$ For fixed¨$\nn, \rr$, the generating functions for
the sequences $\ccc\nn\dd$ and $\hhh\rr\nn\dd$ are
rational.

$(ii)$ For fixed¨$\nn, \rr$, the numbers $\ccc\nn\dd$ and
$\hhh\rr\nn\dd$ admit expressions of the form
\begin{equation} \label{E:General}
P_1(\dd) \ev_1^\dd + \cdots +
P_\ir(\dd) \ev_\ir^\dd.
\end{equation} 
where $\ev_1$, \dots, $\ev_\ir$ are the non-zero
eigenvalues of¨$M_\nn$ and $P_1, \dots,
P_\ir$ are polynomials with $\deg(P_\ind)$ at
most the multiplicity of¨$\ev_\ind$ for¨$M_\nn$.
\end{coro}

As the matrix¨$M_\nn$ is an $\nn! \times \nn!$ matrix,
completing the computation is not so easy, even for small
values of¨$\nn$. Actually, it is shown in¨\cite{Dhi} how
to replace¨$M_\nn$ with a smaller matrix¨$\overline
M_\nn$ of size $p(\nn) \times p(\nn)$, where $p(\nn)$ is
the number of partitions of¨$\nn$. The property is
connected with classical results by Solomon about the
descents of permutations¨\cite{Sol}. With such methods,
one easily obtains the values listed in
Table¨\ref{T:Values}.

\begin{table}[htb]
$$\begin{tabular}{c|r|r|r|r|r|r|r}
\quad$\dd$
& 0 
&1&2&3&4&5&6\\
\hline
\vrule width0pt height12pt depth6pt
$\hhh12\dd$
& 1
& 2 & 3 & 4 & 5 & 6 & 7 \\
\hline
\vrule width0pt height12pt 
$\hhh13\dd$
& 1
& 6 & 19 & 48 & 109 & 234 & 487\\
\vrule width0pt depth6pt
$\hhh23\dd$
& 1
& 3 & 7 & 15 & 31 & 63 & 127\\
\hline
\vrule width0pt height12pt 
$\hhh14\dd$
& 1
&24&211&1\,380&8\,077&45\,252&249\,223
\\
$\hhh24\dd$
& 1
&12&83&492&2\,765&15\,240
& 83\,399 \\
\vrule width0pt depth6pt
$\hhh34\dd$
& 1
&4&15&64&309&1\,600 
& 8\,547
\\
\hline
\vrule width0pt height12pt 
$\hhh15\dd$
& 1
 & 120 & 3\,651 & 79\,140 
& 1\,548\,701 
& 29\,375\,460 
& 551\,997\,751 
\\
$\hhh25\dd$
& 1
 & 60 & 1\,501 & 30\,540 
& 585\,811
& 11\,044\,080 
& 207\,154\,921 
\\
$\hhh35\dd$
& 1
 & 20 & 311 & 5\,260 
& 94\,881
& 1\,755\,360 
& 32\,741\,851 
\\
\vrule width0pt depth6pt
$\hhh45\dd$
& 1
 & 5 & 31 & 325 
& 4\,931 
& 86\,565 
& 1\,590\,231 
\\
\hline
\vrule width0pt height12pt 
$\hhh16\dd$
& 1
& 720 
&\!90\,921 
&\!7\,952\,040 
&\!634\,472\,921
&\!49\,477\,263\,360
&\!3\,836\,712\,177\,121
\\
$\hhh26\dd$
& 1
& 360 
&\!38\,559 
&\!3\,228\,300 
&\!254\,718\,389
&\!19\,808\,530\,620
&\!1\,535\,016\,069\,499 \\
$\hhh36\dd$
& 1
 & 120 & 8\,727 & 649\,260 
& 49\,654\,757
& 3\,831\,626\,580 
& 296\,361\,570\,667\\
$\hhh46\dd$
& 1
 & 30 & 1\,075 & 61\,620 
& 4\,387\,195
& 332\,578\,230 
& 25\,612\,893\,355\\
$\hhh56\dd$
& 1
 & 6 & 63 & 1\,956
& 116\,423
& 8\,448\,606 
& 643\,888\,543
\end{tabular}$$
\bigskip
\caption{\smaller First values of $\hhh\rr\nn\dd$
for $1 \le \rr < \nn$---the value is¨$1$ for $\rr
\ge \nn$. For instance, we read that the number of
$3$-strand braids of degree at most¨$2$, \ie,
$\hhh132$, is¨$19$---as was seen in
Example¨\ref{X:Height}---while the maximal number
of¨$\ss3$'s in a $\s$-positive word drawn in¨$\GG44$,
\ie, $\ccc44$, which is $\hhh344-1$
according to Proposition¨\ref{P:MainHeight}, is¨$308$.}
\label{T:Values}
\end{table}

Using the reduced matrices 
$\overline M_3 = \left(\begin{matrix}1&0&0\\4&2&0\\1&1&1
\end{matrix}\right)$ and 
$\overline M_4 =
\left(\begin{matrix}1&0&0&0&0\\11&4&1&0&0\\5&3&2&1&0\\
6&4&2&2&0\\1&1&1&1&1
\end{matrix}\right)$, we obtain the
following explicit form for¨\eqref{E:General} involving
the non-zero eigenvalues of¨$M_3$, namely $1$ (double),
$2$ and of¨$M_4$, namely $1$ (double), $2$, and $3 \pm
\sqrt6$:

\begin{prop} \label{P:Values34}
Let $\ev_1 = 3 + \sqrt{6}$ and $\ev_2 = 3 - \sqrt{6}$.
Then, for $\dd \ge 1$, we have
\begin{align*}
\hhh13\dd &= 8\cdot 2^\dd - 3 \dd - 7,\\
\hhh23\dd &= \ccc3\dd+1 = 2 \cdot 2^\dd - 1,\\
\hhh14\dd &= 
\frac3{20}(32 + 13\sqrt6) \ev_1^\dd  
+ \frac3{20}(32 - 13\sqrt6) \ev_2^\dd 
- \frac{128}5 \cdot\, 2^\dd
+  6\dd + 17,\\
\hhh24\dd &= 
\frac1{20}(32 + 13\sqrt6) \ev_1^\dd 
+ \frac1{20}(32 - 13\sqrt6) \ev_2^\dd
- \frac{16}5 \cdot\, 2^\dd 
+ 1,\\
\hhh34\dd &= \ccc34 +1 = 
\frac1{20}(4 + \sqrt6) \ev_1^\dd
+ \frac1{20}(4 - \sqrt6) \ev_2^\dd 
+ \frac85 \cdot\, 2^\dd - 1.
\end{align*}
\end{prop}

The main interest of the above formulas are to show that
each of the involved parameters has an exponential growth
with respect to¨$\dd$, in $O(2^\dd)$ for $\nn = 3$, and in
$O((3+\sqrt{6})^\dd)$ for $\nn = 4$. For practical
purposes, it may be more convenient to resort to inductive
formulas, for instance
\begin{gather}
\label{E:Recur3}
\hhh13\dd = 2 \hhh13{\dd-1} + 3 \dd + 1,\\
\label{E:Recur4}
\hhh14\dd = 6\hhh14{\dd-1} - 3\hhh14{\dd-2}
+ 32 \cdot 2^\dd - 12 \dd - 34,
\end{gather}
together with initial values $\hhh130 = \hhh140 = 1$,
$\hhh141 = 24$ (or $\hhh14{-1}= 0$).

\subsection{Small values of¨$\dd$}

Another approach is to keep¨$\dd$ fixed and let¨$\nn$ vary.
Once again, we only mention a few results, and refer the
reader to¨\cite{Dhi} for the proofs and additional
comments. For $\dd = 1$, it is easy to determine all
values:

\begin{prop} [\cite{Dhi}]
For $\nn \ge \rr \ge 1$, we have
\begin{equation*} \label{E:K1}
\hhh\rr\nn{} = \frac{\nn!}{\rr!}.
\end{equation*}
\end{prop}

For $\dd = 2$, it is easier to complete the
computation for $\nn - \rr$ rather than¨$\rr$ fixed. 

\begin{prop} [\cite{Dhi}] \label{P:K2}
For $\nn \ge \rr \ge 1$, we have
\begin{align*} 
\hhh{\nn-\rr}\nn2 
&= \rr!\,(\rr+1)^\nn + \sum_{i = 1}^\rr P_i(\nn) \,
i^{\nn-\rr+i-1},\\
\intertext{for some polynomial $P_i$ of degree at
most¨$\rr - i + 1$. The values for $\rr =
1, 2$ are}
\hhh{\nn-1}\nn2 &= 2^\nn - 1, \\
\hhh{\nn-2}\nn2 &= 2 \cdot 3^\nn - (\nn + 6) \cdot
2^{\nn-1} + 1.
\end{align*}
\end{prop}

For $\rr$ fixed, no general formula is known. Let us
mention the case of¨$\hhh1\nn2$, which follows from
results of¨\cite{CSV}:

\begin{prop} [\cite{Dhi}]
The numbers¨$\hhh1\nn2$ are determined by the
induction
\begin{equation*}
\hhh102=1, \qquad 
\hhh1\nn2 = \sum_{\ind=0}^{\nn-1} (-1)^{\nn + \ind
+1} {\nn \choose \ind}^2 
\hhh1\ind2.
\end{equation*}
Their double exponential generating function is, with 
$J_0(\xx)$ is the Bessel function,
\begin{equation*}
\sum_{\nn = 0}^\infty \hhh1\nn2
\frac{z^\nn}{\nn!^2} = 
\bigg( \sum_{\nn = 0}^\infty (-1)^\nn
\frac{z^\nn}{\nn!^2}
\bigg)\inv = \frac1{J_0(\sqrt{z})}.
\end{equation*}
\end{prop}

Finally, for $\dd = 3$, the computation can be completed
at least in the case $\nn - \rr = 1$:  

\begin{prop} [\cite{Dhi}]
For¨$\nn \ge 1$, we have, with $e = \exp(1)$,
\begin{equation*} \label{E:KKK3}
\hhh{\nn-1}\nn3 = \sum_{\ind=0}^{\nn-1}
\frac{\nn!}{\ind!} = \lfloor\nn! e\rfloor - 1.
\end{equation*}
\end{prop}

Using Proposition¨\ref{P:MainHeight}, we deduce the
following explicit values for $\ccc\nn\dd$, \ie, for the
maximal number of occurrences of¨$\ss{\nn-1}$ is a
$\s$-positive word drawn in the Cayley graph
of¨$\DD\nn\dd$:
$$\ccc\nn{} = \nn - 1, \qquad
\ccc\nn2 = 2^\nn - 2, \qquad
\ccc\nn3 = \sum_{\ind=0}^{\nn-1}
\frac{\nn!}{\ind!} - 1 = \lfloor\nn! e\rfloor - 2.
$$
The formulas listed above show that a number
of different induction schemes appear, suggesting that the
combinatorics of normal sequences of braids is very rich. 

\section{A complete description of $\DDDD3\dd$}
\label{S:N3}

Our ultimate goal would be a complete description
of each chain $\DDDD\nn\dd$, this typically meaning
that we are able to explicitly specify the increasing
enumeration of its elements. This goal remains out of
reach in the general case, but we shall show now how the
process can be completed in the case $\nn = 3$. The
counting formulas of Section¨\ref{S:Deter} play a key
role in the construction, and, in particular, the Pascal
triangle of Table¨\ref{T:Triangle} below is directly
connected with the $2^\dd$¨factor in the inductive
formulas of Proposition¨\ref{P:Values34}.  As an
application, we deduce a new proof of Property¨C and of
the well-ordering property, hence a complete
re-construction of the braid ordering in the case $\nn
= 3$.

The general principle is to make the decomposition of
Corollary¨\ref{C:Structure} explicit. The latter shows
that, for all¨$\nn$ and¨$\dd$, the chain $\DDDD\nn\dd$
can be decomposed into $\ccc\nn\dd$¨subintervals
each of which is a copy of some fragment
of¨$\DDDD{\nn-1}\dd$. Moreover, the approach of
Section¨\ref{S:Deter} suggests an induction on¨$\dd$
as well, so, finally, we are led to looking for a
description of¨$\DDDD\nn\dd$ in terms
of¨$\DDDD{\nn-1}\dd$ and¨$\DDDD\nn{\dd-1}$---\ie, in
the current case, a description of $\DDDD3\dd$ in terms
of
$\DDDD2\dd$ and $\DDDD3{\dd-1}$.

\subsection{The braids¨$\Dd\nn\pp$}

The subsequent construction will appeal to a double
series of braid called¨$\Dd\nn\pp$, and we begin with a
few preliminary properties of these braids.

\begin{defi}
For $\nn \ge 2$, let $\ss{\nn, 1}$ and $\ss{1, \nn}$
respectively denote the braid words
$\ss{\nn-1} \ss{\nn-2} \dots \ss1$ and $\ss1 \ss2
\dots \ss{\nn-1}$. For $\pp \ge 0$, we define 
$\Ddt\nn\pp$ to  be (the braid represented by) the
length¨$\pp$ prefix of the right infinite
word¨$(\ss{\nn, 1}\ss{1, \nn})^\infty$, and $\Dd\nn\pp$ to
be (the braid represented by) the length¨$\pp$ suffix of
the left infinite word¨$^\infty(\ss{\nn, 1}\ss{1, \nn})$. 
\end{defi}

For instance, we find
$\Dd30 = 1$, $\Dd31 = \b$, $\Dd32 = \a\b$, \dots,
$\Dd34 = \b\a\a\b$, \dots, $\Dd37 = \a\a\b\b\a\a\b$, \etc\
Similarly, we have $\Dd46 = \c\b\a\a\b\c$, and, more
generally, $\Dd\nn{2\nn-2} = \Ddt\nn{2\nn-2}= \ss{\nn,
1}\ss{1, \nn}$. Note that, as a word,
$\Dd\nn\pp$ is obtained by reversing the order of the
letters in¨$\Ddt\nn\pp$.

\begin{lemm}
For $\nn \ge 2$ and $\pp, \qq \ge 0$ satisfying $\pp
+ \qq = \dd(\nn-1)$, we have
\begin{equation} \label{E:Deltadelta}
\Dd\nn\pp \, \DD{\nn-1}\dd \, \Ddt\nn\qq = \DD\nn\dd.
\end{equation}
\end{lemm}

\begin{proof}
We first prove using induction on¨$\dd$ the relation
\begin{equation} \label{E:Deltak0}
\Dd\nn{\dd(\nn-1)} \, \DD{\nn-1}\dd = \DD\nn\dd,
\end{equation}
\ie, \eqref{E:Deltadelta} with $\qq =
0$. For $\dd = 0$, \eqref{E:Deltak0} reduces to $1 =
1$. Assume $\dd \ge 1$.  By definition,
$\Dd\nn{\dd(\nn-1)}$ is $\ss{\nn,1} \, 
\Dd\nn{(\dd-1)(\nn-1)}$ for
$\dd$¨odd, and is $\ss{1, \nn} \,
\Dd\nn{(\dd-1)(\nn-1)}$ for
$\dd$¨even, so, in any case, we can write
$$\Dd\nn{\dd(\nn-1)} = \flip\nn^{\dd-1}(\ss{1, \nn}) \,
\Dd\nn{(\dd-1)(\nn-1)},$$ 
where we recall $\flip\nn$ denotes the flip
automorphism of¨$B_\nn$ that exchanges¨$\ss i$ and
$\ss{\nn - i}$. Using the induction hypothesis
and¨\eqref{E:Delta}, we find
\begin{align*}
\Dd\nn{\dd(\nn-1)} \, \DD{\nn-1}\dd 
&= \flip\nn^{\dd-1}(\ss{1, \nn}) \, \Dd\nn{(\dd-1)(\nn-1)} \,
\DD{\nn-1}{\dd-1} \, \D_{\nn-1}\\
&= \flip\nn^{\dd-1}(\ss{1, \nn}) \, \DD n{\dd-1} \, \D_{\nn-1}
= \DD n{\dd-1} \, \ss{1, \nn} \, \D_{\nn-1}
= \DD n{\dd-1} \, \D_\nn 
= \DD\nn\dd.
\end{align*} 

We return to the general case
of¨\eqref{E:Deltadelta}. For $\dd$ even, we have
$\Dd\nn{\dd(\nn-1)} = \Ddt n{\dd(\nn-1)}$, hence
$\Ddt\nn\qq \, \Dd\nn\pp = \Dd\nn{\dd(\nn-1)}$. If
$\dd$ is odd, we have $\Dd\nn{\dd(\nn-1)} =
\flip\nn(\Ddt n{\dd(\nn-1)})$, which implies
$\flip\nn(\Ddt\nn\qq) \, \Dd\nn\pp = \Dd\nn{\dd(\nn-1)}$. So $\flip\nn^\dd(\Ddt\nn\qq) \, \Dd\nn\pp = \Dd\nn{\dd(\nn-1)}$ holds in both
cases. Now, using¨\eqref{E:Deltak0}, we find
\begin{align*}
\flip\nn(\Ddt\nn\qq) \, \Dd\nn\pp \, \DD{\nn-1}\dd \, \Ddt\nn\qq
&= \Dd\nn{\dd(\nn-1)} \, \DD{\nn-1}\dd \, \Ddt\nn\qq
= \DD\nn\dd \, \Ddt\nn\qq
= \flip\nn(\Ddt\nn\qq) \, \DD\nn\dd,
\end{align*} 
from which we deduce¨\eqref{E:Deltadelta} by cancelling 
$\flip\nn(\Ddt\nn\qq)$ on the left.
\end{proof}

\begin{lemm}
For $1 \le i \le \nn-2$ we have
\begin{equation}
\label{E:Commut}
\Dd\nn{\dd(\nn-1)} \, \ss i = \ss{i+e} \,
\Dd\nn{\dd(\nn-1)}
\end{equation}
with $e = 0$ if $\dd$ is even, and $e = 1$ if $\dd$ is odd.
\end{lemm}

\begin{proof}
For $1 \le i \le \nn-2$, we have
\begin{equation}
\ss{1,n} \, \ss i = \ss{i+1} \, \ss{1,n},
\mbox{\quad and \quad}
\ss{\nn, 1} \, \ss{i+1} = \ss i \, \ss{\nn, 1}, 
\end{equation}
as an easy induction shows. This implies
$\ss{\nn,1} \, \ss{1,n} \, \ss i = \ss i \,
\ss{\nn,1} \, \ss{1,n}$, and therefore $(\ss{\nn,1}
\, \ss{1,n})^\dd \, \ss i = \ss i \,
(\ss{\nn,1} \, \ss{1,n})^\dd$, \ie, 
$\Dd\nn{2\dd(\nn-1)} \, \ss i = \ss i \,
\Dd\nn{2\dd(\nn-1)}$, for every¨$\dd$. On the other
hand, we have 
$\Dd\nn{(2\dd+1)(\nn-1)} = \ss{1,n} \, \Dd\nn{2\dd(\nn-1)}$, hence
$$\Dd\nn{(2\dd+1)(\nn-1)} \, \ss i 
= \ss{1,n} \, \ss i \, \Dd\nn{2\dd(\nn-1)} 
= \ss{i+1} \, \ss{1,n} \, \Dd\nn{2\dd(\nn-1)} =
\ss{i+1} \, \Dd\nn{(2\dd+1)(\nn-1)},$$
as was expected.
\end{proof}

\subsection{A Pascal triangle}

We shall now construct for every¨$\dd$ a sequence of
positive braids¨$\SS3\dd$ that will turn out to be
the increasing enumeration of¨$\DDDD3\dd$.
The construction relies on an induction similar to
a Pascal triangle. In order to make it easily
understandable, it is convenient to start with a
construction in the (trivial) cases¨$\nn =
1$ and¨$\nn = 2$.

As $B_1$ is the trivial group, then for every¨$\dd$
there is exactly one element of degree at
most¨$\dd$, namely¨$1$, and we can state:

\begin{prop}
Let $\SS1\dd$ be defined for
$\dd \ge 0$ by
\begin{equation} \label{E:Def1}
\SS1\dd = (1).
\end{equation}
Then $\SS1\dd$ is the increasing enumeration
of¨$\DDD1\dd$.
\end{prop}

The group¨$B_2$ is the rank¨$1$ free group
generated by¨$\ss1$. The fundamental braid¨$\D_2$
is just¨$\ss1$, and the braids of degree at
most¨$\dd$, \ie, the divisors of¨$\DD2\dd$, consist
of the $\dd+1$ braids¨$1, \ss1, \ldots, \sss1\dd$.
On the other hand, we have $\ss{1, 2} = \ss{2, 1} =
\ss1$, and $\Dd1i = \sss1i$ for every¨$i$. 

\begin{nota}
If $S_1, S_2$ are sequences (of braids), we denote
by $S_1 + S_2$ the concatenation of¨$S_1$ and¨$S_2$,
\ie, the sequence obtained by appending¨$S_2$ after¨$S_1$.
If $S$ is a sequence of braids, and $\xx$ is a
braid, we denote by¨$\xx S$ the translated sequence
obtained by multiplying each entry in¨$S$ by¨$\xx$
on the left. 
\end{nota}

With these conventions, the sequence $(1, \ss1,
\dots, \ss1^\dd)$ can be expressed as $\Dd20(1)
+ \Dd21(1) + \cdots + \Dd2\dd(1)$, and we can
state:

\begin{prop} \label{P:Enum3}
Let $\SS2\dd$ be defined for $\dd \ge 0$ by
\begin{equation} \label{E:Def2}
\SS2\dd = \Dd20\SS1\dd + \Dd21\SS1\dd + \cdots + \Dd2\dd\SS1\dd.
\end{equation}
Then $\SS2\dd$ is the increasing enumeration
of¨$\DDD2\dd$.
\end{prop}

We repeat the process for¨$\nn = 3$, introducing
a sequence¨$\SS3\dd$ by a definition similar
to¨\eqref{E:Def2} that involves $\SS2\dd$
and¨$\SS3{\dd-1}$. The result we shall prove is:

\begin{prop} \label{P:Main3}
Let $\SS3\dd$ be defined for $\dd \ge 0$ by
\begin{align} \label{E:Def3}
\SS3\dd =  \Dd30\SS2 \dd + \SS3{\dd, 1}
+ \Dd31\SS2 \dd + \cdots 
+ \Dd3{2\dd-1}\SS2 \dd + \SS3{\dd, 2\dd} + \Dd3{2\dd} \SS2 \dd,
\end{align}
where $\SS3{\dd, 1}, \cdots, \SS3{\dd, 2\dd}$ are 
defined by $\SS3{\dd, 1} =
\SS3{\dd, 2\dd} =
\emptyset$ and, for $2 \le \ppp \le 2\dd-1$, 
\begin{equation*}
\SS3{\dd, \ppp} = 
\begin{cases}
\phantom{\ss2}\ss1(\SS3{\dd-1, \ppp-1} +
\Dd3{\ppp-1}\SS2{\dd-1} +
\SS3{\dd-1, \ppp}) & \mbox{for $\ppp = 0 \pmod 4$,}\\
\ss2\ss1(\SS3{\dd-1, \ppp-2} + \Dd3{\ppp-1}\SS2{\dd-1} +
\SS3{\dd-1, \ppp-1}) & \mbox{for $\ppp = 1 \pmod 4$,}\\
\phantom{\ss2}\ss2(\SS3{\dd-1, \ppp-1} +
\Dd3{\ppp-1}\SS2{\dd-1} +
\SS3{\dd-1, \ppp}) & \mbox{for $\ppp = 2 \pmod 4$,}\\
\ss1\ss2(\SS3{\dd-1, \ppp-2} + \Dd3{\ppp-1}\SS2{\dd-1} +
\SS3{\dd-1, \ppp-1}) & \mbox{for $\ppp = 3 \pmod 4$.}
\end{cases}
\end{equation*}
Then $\SS3\dd$ is the increasing enumeration
of¨$\DDD3\dd$.
\end{prop}

The general scheme is illustrated in 
Table¨\ref{T:Triangle}: the sequence¨$\SS3\dd$ is
constructed by starting with $2\dd+1$¨copies
of¨$\SS2\dd$ translated by¨$\Dd30$, \dots,
$\Dd3{2\dd}$ and inserting (translated copies
of) fragments of the previous sequence¨$\SS3{\dd-1}$.  

\begin{table}[t]
\begin{picture}(130,47)(1, 0)
\put(0,10){$\Dd30 \SS23$}
\put(10,10){$(\SS3{3,1})$}
\put(20,10){$\Dd31 \SS23$}
\put(32,10){$\SS3{3,2}$}
\put(40,10){$\Dd32 \SS23$}
\put(52,10){$\SS3{3,3}$}
\put(60,10){$\Dd33 \SS23$}
\put(72,10){$\SS3{3,4}$}
\put(80,10){$\Dd34 \SS23$}
\put(92,10){$\SS3{3,5}$}
\put(100,10){$\Dd34 \SS23$}
\put(110,10){$(\SS3{3,6})$}
\put(120,10){$\Dd36 \SS23$}
\put(10,9){$\underbrace{\hbox to 28mm{\hfill}}$}
\put(14,4){$\ss2\cdot$}
\put(29,4){$\ss1\ss2\cdot$}
\put(17.5,3){$\swarrow$}
\put(26.5,3){$\searrow$}
\put(52,9){$\underbrace{\hbox to 26mm{\hfill}}$}
\put(56,4){$\ss1\cdot$}
\put(70,4){$\ss2\ss1\cdot$}
\put(59,3){$\swarrow$}
\put(68,3){$\searrow$}
\put(92,9){$\underbrace{\hbox to 27mm{\hfill}}$}
\put(96,4){$\ss2\cdot$}
\put(111,4){$\ss1\ss2\cdot$}
\put(99,3){$\swarrow$}
\put(108,3){$\searrow$}
\put(14,0){$\cdots$}
\put(31,0){$\cdots$}
\put(54,0){$\cdots$}
\put(71,0){$\cdots$}
\put(95,0){$\cdots$}
\put(112,0){$\cdots$}

\put(20,22){$\Dd30 \SS22$}
\put(30,22){$(\SS3{2,1})$}
\put(40,22){$\Dd31 \SS22$}
\put(52,22){$\SS3{2,2}$}
\put(60,22){$\Dd32 \SS22$}
\put(72,22){$\SS3{2,3}$}
\put(80,22){$\Dd33 \SS22$}
\put(90,22){$(\SS3{2,4})$}
\put(100,22){$\Dd34 \SS22$}
\put(30,21){$\underbrace{\hbox to
28mm{\hfill}}$}
\put(36,16){$\ss2\cdot$}
\put(50,16){$\ss1\ss2\cdot$}
\put(38,14.5){$\swarrow$}
\put(47,14.5){$\searrow$}
\put(72,21){$\underbrace{\hbox to
27mm{\hfill}}$}
\put(77,16){$\ss1\cdot$}
\put(91,16){$\ss2\ss1\cdot$}
\put(79,14.5){$\swarrow$}
\put(88,14.5){$\searrow$}

\put(40,34){$\Dd30 \SS21$}
\put(50,34){$(\SS3{1,1})$}
\put(60,34){$\Dd31 \SS21$}
\put(70,34){$(\SS3{1,2})$}
\put(80,34){$\Dd32 \SS21$}
\put(50,33){$\underbrace{\hbox to
28mm{\hfill}}$}
\put(56,28){$\ss2\cdot$}
\put(70,28){$\ss1\ss2\cdot$}
\put(58,26.5){$\swarrow$}
\put(67,26.5){$\searrow$}

\put(60,44){$\Dd30 \SS20$}
\end{picture}
\smallskip
\caption{\smaller The inductive construction
of¨$\SS3\dd$ as a Pascal triangle: the subsequence¨$
\SS3{\dd,\ppp}$ is obtained by (translating and)
concatenating the previous subsequences¨$
\SS3{\dd-1,\ppp-1}$ and¨$
\SS3{\dd-1,\ppp}$, or¨$\SS3{\dd-1,\ppp-2}$ and¨$
\SS3{\dd-1,\ppp-1}$, depending on the parity of¨$\ppp$;
the parenthesized sequences are empty; if we forget about
the subsequences¨$\Dd3\qq\SS2\dd$, we have the
Pascal triangle.}
\label{T:Triangle}
\end{table}

\begin{exam} \label{X:SS3}
The difference between the definition of¨$\SS3\dd$
in¨\eqref{E:Def3} and that of¨$\SS2\dd$ in¨\eqref{E:Def2} is
the insertion of the additional factors $ \SS3{\dd,\ppp}$
between the consecutive terms¨$\Dd3\qq\SS2\dd$. Because
$\SS3{\dd, 1}$ and $\SS3{\dd, 2\dd}$ are empty, the difference
occurs for $\dd \ge 2$ only. The first values are:

$\SS30 = \Dd30 \SS20 = (1)$,

$\SS31 = \Dd30 \SS21 + \SS3{1,1} + \Dd31 \SS21 +
\SS3{1,2} + \Dd32 \SS21$

\hspace{1cm} $= (1, \a) + \emptyset
+ \b(1, \a) + \emptyset 
+ \a\b(1, \a)$

\hspace{1cm} $= (1, \a, \b, \b\a, \a\b, \a\b\a)$,

$\SS32  = \Dd30 \SS22 + \SS3{2,1} +
\Dd31 \SS22 + \SS3{2,2} +
\Dd32 \SS22 + \SS3{2,3} +
\Dd33 \SS22 + \SS3{2,4} +
\Dd34 \SS22$

\hspace{1cm} $ = (1, \a, \a\a) + \emptyset
+ \b(1, \a, \a\a) + \b(\b, \b\a)
+ \a\b(1, \a, \a\a)$

\hspace{4cm} $ + \a\b(\b, \b\a)
+ \a\a\b(1, \a, \a\a) + \emptyset 
+ \b\a\a\b(1, \a, \a\a)$

\hspace{1cm} $ = (1, \a, \a\a, \b, \b\a, \b\a\a,
\b\b, \b\b\a,
\a\b, \a\b\a, \a\b\a\a,
\a\b\b, \a\b\b\a, \a\a\b$, 

 \hspace{6cm} $\a\a\b\a, \a\a\b\a\a, \b\a\a\b,
\b\a\a\b\a, \b\a\a\b\a\a)$.\\
It is easy to check directly that the sequence
$\SS3\dd$ provides the increasing enumeration
of¨$\DDD3\dd$ for $\dd = 0, 1, 2$. 
\end{exam}

The proof of Proposition¨\ref{P:Main3} will be
split into several pieces, each of which is
established using an induction on the degree¨$\dd$.

\begin{lemm} \label{L:Inclusion3}
All entries in¨$\SS3\dd$ are divisors of¨$\DD3\dd$.
\end{lemm}

\begin{proof}
The result is true for $\dd = 0$. Assume $\dd \ge 1$. By
construction, each entry in¨$\SS3\dd$ either is of the form
$\Dd3\qq \ss1^\expo$ with $0 \le \qq \le 2\dd$ and $0 \le
\expo \le \dd$, or belongs to some
subsequence¨$\SS3{\dd,\ppp}$ with $2 \le \ppp \le 2\dd -
1$. In the first case,
$\Dd3\qq\ss1^\expo$ is a right divisor of $\Dd3{2\dd}
\ss1^\expo$, which itself is a left divisor of
$\Dd3{2\dd} \ss1^\dd$. By¨\eqref{E:Deltadelta}, the latter
is¨$\DD3\dd$. Hence each $\Dd3\qq
\ss1^\expo$ is a divisor of¨$\DD3\dd$. As for the
entries coming from some subsequence¨$\SS3{\dd,\ppp}$,
by definition they are of the form¨$\xx\yy$ with $\xx$
one of $\ss2, \ss1\ss2, \ss1,
\ss2\ss1$ and $\yy$ an entry in¨$\SS3{\dd-1}$. Then
$\xx$ is a divisor of¨$\D_3$, while, by induction
hypothesis, $\yy$ is a divisor of¨$\DD3{\dd-1}$, so
$\xx\yy$ is a divisor of¨$\DD3\dd$.
\end{proof}

\begin{lemm} \label{L:Card3}
The length of the sequence¨$\SS3\dd$ equals the
cardinality of¨$\DDD3\dd$.
\end{lemm}

\begin{proof}
Let $\ell_\dd$ denote the length of¨$\SS3\dd$.
Computing¨$\ell_\dd$ is not very difficult. However,
there is need to do it, which amounts to solve a
recursive formula unnecessarily. Indeed, we saw in
Section¨\ref{S:Deter} that the cardinality¨$\hhh13\dd$
of¨$\DDD3\dd$ obeys the inductive rule¨\eqref{E:Recur3}.
So it will be enough to check that
$\ell_\dd$ satisfies the relation
\begin{equation} \label{E:Double}
\ell_\dd = 2 \ell_{\dd-1} + 3\dd + 1,
\end{equation}
and starts from the initial¨$\ell_1 = 6$ (or $\ell_0 =
1$). The latter point was checked in Example¨\ref{X:SS3}.

Now Table¨\ref{T:Triangle} shows that most entries
in¨$\SS3{\dd-1}$ give rise to two entries in¨$\SS3\dd$.
More precisely, each entry of¨$\SS3{\dd-1}$ not
belonging to a factor of the form¨$\Dd3{2\qq}
\SS2{\dd-1}$ gives rise to two entries in¨$\SS3\dd$,
and, conversely, each entry in¨$\SS3\dd$ not belonging
to a factor¨$\Dd3{\qq}
\SS2{\dd}$ comes from such an entry in¨$\SS3{\dd-1}$.
There are $\dd$¨factors¨$\Dd3{2\qq} \SS2{\dd-1}$
in¨$\SS3{\dd-1}$, each of length¨$\dd$, and
$2\dd+1$¨factors¨$\Dd3{2\qq} \SS2{\dd}$ in¨$\SS3\dd$,
each of length¨$\dd+1$. So we obtain
$$\ell_\dd - (2\dd+1)(\dd+1) = 2(\ell_{\dd-1} - \dd^2),$$
which gives¨\eqref{E:Double}.
\end{proof}

At this point, we cannot (yet) conclude that each
divisor of¨$\DD3\dd$ occurs exactly once in¨$\SS3\dd$,
as there could be some repetitions.

\subsection{A quotient-sequence for¨$\SS3\dd$}

Our next aim is to show that $\SS3\dd$ is
$\sm$-increasing. To this end, we shall explicitly
determine the quotient of adjacent entries
in¨$\SS3\dd$, \ie, we shall specify a quotient-sequence
for¨$\SS3\dd$ in the sense of Definition¨\ref{D:Witness}.

A preliminary step consists in determining the first and
the last entries of the sequence¨$ \SS3{\dd,\ppp}$.
For¨$S$ a nonempty sequence, we denote by $\first S$
(\resp $\last S$) the first (\resp last) entry in¨$S$.

\begin{lemm} \label{L:FirstLast}
For $1 < i < 2\dd$, we have
\begin{equation}
\first{ \SS3{\dd,\ppp}} = \Dd3{\ppp-1} \, \ss2, 
\mbox{\quad and \quad}
\last{ \SS3{\dd,\ppp}} \, \ss2 = \Dd3\ppp \, \sss1\dd. 
\end{equation}
\end{lemm}

\begin{proof}
The result is vacuously true for $\dd = 0, 1$. Assume $\dd \ge
2$ with $\ppp  = 0 \pmod 4$. Using the definition, the
induction hypothesis, and¨\eqref{E:Commut}, we find
\begin{gather*}
\first{ \SS3{\dd,\ppp}} 
= \ss1 \, \first{ \SS3{\dd-1,\ppp-1}} = 
\ss1 \, \Dd3{\ppp-2} \, \ss2
= \Dd3{\ppp-1} \, \ss2, \\
\last{ \SS3{\dd,\ppp}} \, \ss2 
= \ss1 \, \last{ \SS3{\dd-1,\ppp}} \, \ss2 
=  \ss1 \, \Dd3\ppp \, \sss1{\dd-1} 
= \Dd3\ppp \, \sss1\dd.
\end{gather*}
Similarly, for $\ppp  = 1 \pmod 4$, we have
\begin{gather*}
\first{ \SS3{\dd,\ppp}} 
= \ss2\ss1 \, \first{ \SS3{\dd-1,\ppp-2}} = 
\ss2\ss1 \, \Dd3{\ppp-3} \, \ss2
= \Dd3{\ppp-1} \, \ss2, \\
\last{ \SS3{\dd,\ppp}} \, \ss2 
= \ss2\ss1 \, \last{ \SS3{\dd-1,\ppp-1}} \, \ss2 
=  \ss2\ss1 \, \Dd3{\ppp-1} \, \sss1{\dd-1} 
= \ss2\, \Dd3{\ppp-1} \, \sss1\dd
= \Dd3\ppp \, \sss1\dd.
\end{gather*}
Then, for $\ppp  = 2 \pmod 4$, we have
\begin{gather*}
\first{ \SS3{\dd,\ppp}} 
= \ss2 \, \first{ \SS3{\dd-1,\ppp-1}} = 
\ss2 \, \Dd3{\ppp-2} \, \ss2
= \Dd3{\ppp-1} \, \ss2, \\
\last{ \SS3{\dd,\ppp}} \, \ss2 
= \ss2 \, \last{ \SS3{\dd-1,\ppp}} \, \ss2 
=  \ss2 \, \Dd3\ppp \, \sss1{\dd-1} 
= \Dd3\ppp \, \sss1\dd.
\end{gather*}
Finally, for $\ppp  = 3 \pmod 4$, we find
\begin{gather*}
\first{ \SS3{\dd,\ppp}} 
= \ss1\ss2 \, \first{ \SS3{\dd-1,\ppp-2}} = 
\ss1\ss2 \, \Dd3{\ppp-3} \, \ss2
= \Dd3{\ppp-1} \, \ss2, \\
\last{ \SS3{\dd,\ppp}} \, \ss2 
= \ss1\ss2 \, \last{ \SS3{\dd-1,\ppp-1}} \, \ss2 
=  \ss1\ss2 \, \Dd3{\ppp-1} \, \sss1{\dd-1} 
=  \ss1\ss2\ss1\ss2 \, \Dd3{\ppp-3} \, \sss1{\dd-1}\\
\hspace{4cm}
=  \ss1\ss1\ss2\ss1 \, \Dd3{\ppp-3} \, \sss1{\dd-1} 
=  \ss1\ss1\ss2 \, \Dd3{\ppp-3} \, \sss1\dd
= \Dd3\ppp \, \sss1\dd.
\end{gather*}
This completes the argument.
\end{proof}

We shall now construct an explicit quotient-sequence
for¨$\SS3\dd$, \ie, a sequence of braid words
representing the quotients of the consecutive entries
of¨$\SS3\dd$. Before doing it for¨$\SS3\dd$, let us
consider the (trivial) cases of¨$\SS1\dd$ and¨$\SS2\dd$.
As $\SS1\dd$ consists of one single entry, it vacuously
admits the empty sequence as a quotient-sequence. As
for¨$\SS2\dd$, we can state:

\begin{lemm}
For $\dd \ge 0$, let $\WW1\dd$ be the empty sequence, and
let¨$\WW2\dd$ be defined by
\begin{equation} \label{E:DDef2}
\WW2\dd = \WW1\dd + (\ss1) + \WW1\dd + \cdots + \WW1\dd
+ (\ss1) + \WW1\dd,
\end{equation}
$\dd$¨times¨$(\ss1)$. Then $\WW2\dd$ is a
quotient-sequence for¨$\SS2\dd$.
\end{lemm}

On a similar way, we shall prove:

\begin{prop} \label{P:Quotient3}
Let $\WW3\dd$ be the sequence defined by
$\WW30 = \emptyset$ and 
\begin{align} \label{E:DDef3}
\WW3\dd =  \WW2\dd 
\!+\! (\ss1^{-\dd}\ss2) 
&+ \WW2\dd 
\!+\! (\ss1^{-\dd}\ss2) 
\!+\! \WW3{\dd, 2} \!+\!
(\ss2\ss1^{-\dd}) \\ 
\notag &+ \WW2\dd 
\!+\! (\ss1^{-\dd}\ss2) 
\!+\! \WW3{\dd, 3} 
\!+\!
(\ss2\ss1^{-\dd}) 
\!+\! \cdots \\  
\notag &+ \WW2\dd 
\!+\! (\ss1^{-\dd}\ss2) 
\!+\! \WW3{\dd, 2\dd-1} \!+\!
(\ss2\ss1^{-\dd}) 
\!+\! \WW2\dd 
\!+\! (\ss2\ss1^{-\dd})
\!+\! \WW2\dd,
\end{align}
\begin{tabbing} 
with \= $\WW3{\dd, 2} = \WW3{\dd, 3} 
= \WW2{\dd-1} 
\!+\! (\ss2\sss1{-\dd+1}) 
\!+\! \WW3{\dd-1, 2},$\\ 
\> \vline width0pt height15pt
$\WW3{\dd, 2\dd-2} =\WW3{\dd, 2\dd-1}
= \WW3{\dd-1, 2\dd-3} 
\!+\! (\sss1{-\dd+1}\ss2) 
\!+\! \WW2{\dd-1}$,\\
and
\> \vline width0pt height15pt
$\WW3{\dd, 2\pppp}  = \WW3{\dd, 2\pppp+1} 
= \WW3{\dd-1,2\pppp-1} 
\!+\! (\sss1{-\dd+1}\ss2) 
\!+\! \WW2{\dd-1} 
\!+\! (\ss2\sss1{-\dd+1}) 
\!+\! \WW3{\dd-1, 2\pppp}$
\end{tabbing}
for $4 \le 2\pppp \le 2\dd-4$. Then
$\WW3\dd$ is a quotient-sequence for¨$\SS3\dd$.
\end{prop}

\begin{exam} \label{X:Quotient}
We find
$\WW31 
 = \WW21 + (\A\b) + \WW21 + (\b\A) + \WW21
 = (\a, \A\b, \a, \b\A, \a)$, and
\begin{align*}
\WW32 
= \WW22 + (\A\A\b) 
&+ \WW22 + (\A\A\b) 
+ \WW3{2, 2} + (\b\A\A)\\
&+ \WW22 +
(\A\A\b) + \WW3{2, 3} + (\b\A\A) + \WW22 +
(\b\A\A) + \WW22
\end{align*}
 with $\WW3{2, 2} = \WW3{2, 3} = \WW21 = (\a)$,
whence
$$\WW32 
 = (\a, \a, \A\A\b, \a, \a, \A\A\b, \a, \b\A\A, \a, \a,
\A\A\b, \a, \b\A\A, \a, \a, \b\A\A, \a, \a).$$
\end{exam}

\begin{proof} [Proof of Proposition¨\ref{P:Quotient3}]
We prove using induction on¨$\dd$ that $\WW3\dd$ is a
quotient-sequence for¨$\SS3\dd$ with the $4\dd-2$¨terms
in¨\eqref{E:DDef3} corresponding to the $4\dd-1$¨nonempty
terms in¨\eqref{E:Def3}---so, in particular, for $2
\le \ppp \le 2\dd-1$, the subsequence¨$\WW3{\dd,\ppp}$ is
a quotient-sequence for ¨$ \SS3{\dd,\ppp}$. The result is
vacuously true for
$\dd=0$.  
Assume $\dd \ge 1$. By definition, the sequence¨$\SS3\dd$
consists of the concatenation of the $2\dd+1$ sequences
$\Dd30 \SS2 \dd, \cdots, \Dd3{2\dd} \SS2 \dd$, in which the $2\dd-2$
sequences $\SS3{\dd, 2}, \dots, \SS3{\dd, 2\dd-1}$ are inserted.
We shall consider these subsequences separately, and then
consider the transitions between consecutive subsequences.

First, $\WW2\dd$ is a quotient-sequence for¨$\SS2\dd$,
hence it is a quotient-sequence for every sequence
$\Dd3\qq
\SS2\dd$ as well, since, by definition, the quotients we
consider are invariant under left translation. Then, by
construction, each subsequence
$\SS3{\dd, 2\pppp}$ or $\SS3{\dd, 2\pppp+1}$ appearing in¨$\SS3\dd$ is
obtained by translating some subsequence¨$S$
of¨$\SS3{\dd-1}$, namely
$$S =  \SS3{\dd-1, 2\pppp-1} + \Dd3{\qq-1} \SS2{\dd-1} +
\SS3{\dd-1, 2\pppp}.$$
By induction hypothesis, the sequence
$$\WW3{\dd-1, 2\pppp-1} +
(\sss1{-\dd+1}\ss2) + \WW2{\dd-1} +
(\ss2\sss1{-\dd+1}) +
\WW3{\dd-1, 2\pppp},$$
which is precisely $\WW3{\dd, 2\pppp}$ and $\WW3{\dd, 2\pppp+1}$
by definition, is a quotient-sequence for¨$S$. The
property remains true in the special cases $\pppp = 1$
and $\pppp = \dd$, which correspond to removing the
initial term¨$\SS3{\dd-1,2\pppp-1}$, and/or the final
term¨$\SS3{\dd-1, 2\pppp}$, respectively. Then
$\WW3{\dd, 2\pppp}$ and $\WW3{\dd, 2\pppp+1}$ are also
quotient-sequences for any sequence obtained from¨$S$ by
a left translation, in particular for¨$\SS3{\dd,
2\pppp}$ and¨$\SS3{\dd, 2\pppp+1}$.

So it only remains to study the transitions between the
consecutive terms in the expression¨\eqref{E:Def3}
of¨$\SS3\dd$, \ie, to compare the last entry in each term
with the first entry in the next term. Four cases are to be
considered, namely the special case of the first two terms
and of the final two terms, and the generic cases of the
transitions from¨$\Dd3\qq \SS2\dd$ to¨$\SS3{\dd,\ppp+1}$
and from¨$\SS3{\dd,\ppp}$ to¨$\Dd3\qq \SS2\dd$.

As for the first two terms, namely $\Dd30 \SS2\dd$ and
$\Dd31 \SS2\dd$, \ie, $\SS2\dd$ and $\ss2\SS2\dd$, the last
entry in¨$\SS2\dd$ is¨$\sss1\dd$, while the first entry
in¨$\ss2\SS2\dd$ is¨$\ss2$, so $\sss1{-\dd}\ss2$ is a quotient.
As for the last two terms, namely $\Dd3{2\dd-1} \SS2\dd$ and
$\Dd3{2\dd} \SS2\dd$, the last entry in¨$\Dd3{2\dd-1}\SS2\dd$
is¨$\Dd3{2\dd-1}\,\sss1\dd$, while the first entry
in¨$\Dd3{2\dd}\SS2\dd$ is¨$\Dd3{2\dd}$. Now,
by¨\eqref{E:Deltadelta}, we have 
$\Dd3{2\dd-1}\,\sss1\dd\ss2 =
\Dd3{2\dd}
\,\sss1\dd$, so $\ss2\sss1{-\dd}$ is an expression of the
quotient.

Consider now the transition from¨$\Dd3\qq \SS2\dd$
to¨$\SS3{\dd,\qq+1}$. The last entry in¨$\Dd3\qq \SS2\dd$
is¨$\Dd3\qq\, \sss1\dd$, while, by Lemma¨\ref{L:FirstLast}, the
first entry in¨$\SS3{\dd,\qq+1}$ is¨$\Dd3\qq \,\ss2$. Hence
$\sss1{-\dd}\ss2$ represents the quotient. Finally,
consider the transition from¨$\SS3{\dd,\ppp}$ to¨$\Dd3\qq
\SS2\dd$. By Lemma¨\ref{L:FirstLast} again, the last
entry¨$\xx$ in¨$\Dd3\qq \SS2\dd$ satisfies¨$\xx\, \ss2
= \Dd3\qq \, \sss1\dd$, while the first entry in¨$\Dd3\qq
\SS2\dd$ is¨$\Dd3\qq$. Hence $\ss2\sss1{-\dd}$
represents the quotient.
\end{proof}

\begin{coro} \label{C:Distinct3}
For each¨$\dd$ the sequence¨$\SS3\dd$
is $\sm$-increasing; so, in particular, it consists of
pairwise distinct braids.
\end{coro}

\begin{proof}
By definition, every word in¨$\WW3\dd$ is
$\s$-positive, hence, by Property¨A, it does not
represent¨$1$.
\end{proof}

As $\SS3\dd$ consists of pairwise distinct divisors
of¨$\DD3\dd$, Lemma¨\ref{L:Card3} implies that every
divisor of¨$\DD3\dd$ occurs exactly once in¨$\SS3\dd$.
Then, as $\SS3\dd$ is $\sm$-increasing, it must be the
increasing enumeration of¨$\DDD3\dd$, and the proof of
Proposition¨\ref{P:Main3} is complete.

\begin{rema}
Once we know that $\SS3\dd$ is the increasing
enumeration of¨$\DDD3\dd$ and that $\WW3\dd$ is a
$\s$-positive quotient-sequence for¨$\SS3\dd$, we can
count the $2$-jumps in¨$\SS3\dd$ and obtain the value
of¨$\hhh23\dd$ directly: this amounts to forgetting
about all¨$\ss1^{\pm1}$ in the construction
of¨$\WW3\dd$, and it is then fairly obvious that there
only remains $2^\dd-2$ times¨$\ss2$.
\end{rema}

\subsection{Larger values of¨$\nn$}

The same construction can be developped for¨$\nn = 4$
and further. The general scheme is clear, namely to
define¨$\SS4\dd$ using an inductive rule
\begin{align} \label{E:Def4}
\SS4\dd =  \Dd40 \SS3\dd + \SS4{\dd, 1}
+ \Dd41 \SS3\dd + \cdots 
+ \Dd4{3\dd-1} \SS3\dd + \SS4{\dd, 3\dd} + \Dd4{3\dd}
\SS3\dd,
\end{align}
where the intermediate factor¨$\SS4{\dd,\ppp}$ is
constructed by concatenating and translating convenient
fragments of¨$\SS4{\dd-1}$. Owing to the inductive
rule¨\eqref{E:Recur4} satisfied by the number of
elements¨$\hhh14\dd$ of¨$\DDD4\dd$, we can expect that the
generic entry of¨$\SS4{\dd-1}$ has to be repeated 6¨times
in¨$\SS4\dd$, but with some entries from¨$\SS4{\dd-2}$
repeated 3¨times only. Once the inductive definition
of¨$\SS4\dd$ is made complete, showing that the sequence
is $\sm$-increasing and counting its entries should be
easy. As we have no complete description so far, we
leave the question open here.

\subsection{A new construction for the linear ordering
of¨$B_3$}

The main interest of the approach described above is not
only to connect the Garside structure of¨$B_\nn$ with
its linear ordering, but also to provide a new
independent construction of the braid ordering, at least
in the case of¨$B_3$ as the latter is the only one for
which the construction was completed so far.

As was recalled in the introduction, the existence of
the linear ordering of braids relies on two properties
of braids, namely Property¨A and Property¨C. These
properties have received a number of independent
proofs¨\cite{Dgr}. In particular, Property¨A has now
a very short proof based on Dynnikov's coordinization
for triangulations of a punctured disk (\cite{Dgr},
Chapter¨9). As for Property¨C, no really simple proof
exists so far: not to mention the initial argument
involving self-distributive algebra, the combinatorial
proofs based on handle reduction or on Burckel's
uniform tree approach, as well as the geometric proofs
based on standardization of curve diagrams require some
care. So, at the moment, one can estimate that the
optimal proof of Property¨C is still missing.

A direct application of our construction of the
sequence¨$\SS3\dd$ is

\begin{prop}
Property¨C holds for¨$B_3$, \ie, every non-trivial
$3$-braid admits a $\s$-positive or a $\s$-negative
expression.
\end{prop}

\begin{proof}[New proof]
We take as an hypothesis that Property¨A is true, so
that the relation¨$\sm$ is a partial ordering, but we do
not assume that $\sm$ is linear. As every braid in¨$B_3$
is the quotient of two positive braids in¨$B_3^+$,
proving Property¨C for¨$B_3$ amounts to proving that, if
$\bx, \by$ are arbitrary elements of¨$B_3^+$, then the
quotient $\bx\inv \by$ admits a $\s$-positive or a
$\s$-negative expression.

Now the construction of¨$\SS3\dd$ is self-contained, and
so is that of¨$\WW3\dd$. Then, by construction, every
word in¨$\WW3\dd$ is $\s$-positive. As any concatenation
of $\s$-positive words is $\s$-positive, it follows
that, if $\bx, \by$ are any braids occurring
in¨$\bigcup_\dd\SS3\dd$, then the quotient¨$\bx\inv \by$
admits a $\s$-positive or a $\s$-negative expression,
according to whether $\bx$ occurs before or after¨$\by$
in¨$\SS3\dd$. So, in order to conclude that Property¨C is
true, it just remains to check that each positive
$3$-braid occurs in¨$\bigcup_\dd\SS3\dd$. As
every entry of¨$\SS3\dd$ belongs to¨$\DDD3\dd$, this
is equivalent to proving that each divisor of¨$\DD3\dd$
occurs in¨$\SS3\dd$. Property¨A guarantees that the
entries of¨$\SS3\dd$ are pairwise distinct
(Corollary¨\ref{C:Distinct3}), so it suffices to compare
the length of¨$\SS3\dd$ with the cardinality
of¨$\DDD3\dd$, and this is what we
made in Lemma¨\ref{L:Card3}.
\end{proof}

Actually the construction of¨$\SS3\dd$ gives more.
The approach developped by S.\,Burc\-kel in¨\cite{Bus}
consists in introducing a convenient notion of normal
braid words such that every positive braid admits
exactly one normal expression. In the case of $3$¨
strand braids, the definition is as follows. Every
positive $3$¨strand braid word¨$\ww$ can be written as
an alternated product of blocks¨$\ss1^{\expo}$
and¨$\ss2^{\expo}$. Then we define the {\it code}
of¨$\ww$ to be the sequence made by the sizes of these
blocks. To avoid ambiguity, we decide that the last
block is considered to be a block of¨$\ss1$'s, \ie, we
decide that the code of¨$\ss1$ is¨$(1)$, while the code
of¨$\ss2$ is¨$(1, 0)$. For instance, the code
of¨$\ss2^2
\ss1^3\ss2^5$ is $(2, 3, 5, 0)$.

\begin{defi}
A positive $3$¨strand braid word¨$\ww$ is said to be
{\it normal in the sense of Burckel} if its code has the
form $(e_1, \dots, e_\ell)$ with $e_\kk \ge 2$ for $2
\le \kk\le \ell-2$.
\end{defi}

Burckel shows in¨\cite{Bus} that every positive
$3$-braid admits a unique normal expression and,
moreover, that $\bx \sm \by$ holds if and only if the
normal form of¨$\bx$ is $\mathtt{ShortLex}$-smaller than
the normal form of¨$\by$, where $\mathtt{ShortLex}$
refers to the variant of the lexicographic ordering of
sequences in which the length is given priority: $(e_1,
\dots, e_\ell) <_{\mathtt{ShortLex}} (e'_1, \dots,
e'_{\ell'})$ holds for $\ell < \ell'$, and for $\ell =
\ell'$ and $(e_1, \dots, e_\ell)$ lexicographically
smaller than¨$(e'_1, \dots, e'_{\ell'})$. Burckel's
method consists in defining an iterative reduction
process on non-normal braid words. Our current
approach allows for an alternative, simpler method.
First, a direct inspection shows:

\begin{lemm}
Let $\SSS3\dd$ be the sequence of braid words defined by
the inductive rule¨\eqref{E:Def3}. Then $\SSS3\dd$
consists of words that are normal in the sense of
Burckel.
\end{lemm}

Then, by construction, every braid in¨$\SS3\dd$ is
represented by a word of¨$\SSS3\dd$, so, as every
positive $3$-braid occurs in¨$\bigcup\SS3\dd$, we
immediately deduce:

\begin{prop}
Every positive $3$-braid admits an expression that is
normal in the sense of Burckel.
\end{prop}

This in turn enables us to obtain a simple proof for
the following deep, and so far not very well understood
result originally due to Laver¨\cite{Lve} (and
to Burckel \cite{Bus} for the ordinal type):

\begin{coro}
The restriction of¨$\sm$ to¨$B_3^+$ is a well-ordering
of ordinal type¨$\omega^\omega$.
\end{coro}

\begin{proof}
The $\mathtt{ShortLex}$ ordering of sequences of
nonnegative integers is a well-ordering of ordinal
type¨$\omega^\omega$, so its restriction to codes of
normal words in the sense of Burckel is a well-ordering
as well. The type of the latter cannot be less
than¨$\omega^\omega$ as one can easily exhibit an
increasing sequence of length¨$\omega^\omega$.
\end{proof}

Burckel's approach extends to all braid
monoids¨$B_\nn^+$, at the expense of introducing a
convenient notion of normal word and defining an
associated reduction process which is very intricate.
Completing the construction of the sequences¨$\SS4\dd$
and, more generally, $\SS\nn\dd$ along the lines
described above would hopefully allow for an alternative
much simpler approach. In particular, once the correct
definition is given, all subsequent proofs should reduce
to easy inductive verifications.

\end{document}